\documentclass[a4paper,11pt]{article}
\usepackage{microtype}
\usepackage[english,italian]{babel}
\usepackage[dvips]{graphicx}
\usepackage{psfrag}
\usepackage{float}
\usepackage{amsfonts}
\usepackage{amssymb}
\usepackage{amsmath}
\usepackage{amsthm}
\usepackage{enumerate}
\usepackage{MnSymbol}
\usepackage{subfigure}
\usepackage[toc,page]{appendix}

\newcommand{\G}{\Gamma}

\title{Decompositions of complete multipartite\\ graphs via generalized graceful labelings}

\author{Anna Benini and Anita Pasotti \thanks{Dipartimento di Matematica,
Facolt\`a di Ingegneria, Universit\`a degli Studi di Brescia,
Via Valotti, 9, I-25133 Brescia, Italy. email: anna.benini@ing.unibs.it,\ anita.pasotti@ing.unibs.it}}

\date{}
\newtheorem{defi}{Definition}[section]

\newtheorem{rem}[defi]{Remark}
\newtheorem{ex}[defi]{Example}
\newtheorem{thm}[defi]{Theorem}
\newtheorem{cor}[defi]{Corollary}

\begin{document}
\selectlanguage{english}
\maketitle
\selectlanguage{english}

\begin{abstract}
We prove the existence of infinite classes of cyclic $\G$-decompositions
of the complete multipartite graph, $\G$ being a caterpillar,
a hairy cycle or a cycle.
All the results are obtained by the construction of $d$-divisible
$\alpha$-labelings of $\G$, introduced in \cite{APArs}
as a generalization of classical $\alpha$-labelings,
whose existence implies that one of graph decompositions.
\end{abstract}


\section{Introduction}
As usual, we denote by $K_v$ and $K_{m\times n}$ the \emph{complete
graph on $v$ vertices} and the \emph{complete $m$-partite graph with parts of size $n$},
respectively.
For any graph $\G$ we write $V(\G)$ for the set of its vertices and $E(\G)$
for the set of its edges. If $|E(\G)|=e$, we say that $\G$ has \emph{size} $e$.

Given a subgraph $\G$ of a graph $K$, a $\G$-\emph{decomposition of} $K$
is a set of graphs isomorphic to $\G$, called \emph{blocks},  whose edges
partition the edge-set of $K$. Such a decomposition is said to be \emph{cyclic} when it is invariant
under a cyclic permutation of all the vertices of $K$.
A $\G$-decomposition of $K_v$ is also called a $\G$-\emph{system of order} $v$.
For a survey on graph decompositions see \cite{BE}.

The concept of a \emph{graceful labeling} of $\G$, introduced by A. Rosa \cite{R},
is a useful tool working indirectly to determine infinite classes of $\G$-systems.
Rosa proved that
if a graph $\G$ of size $e$ admits a graceful labeling then there exists a cyclic
$\G$-system of order $2e+1$ and if $\G$
admits an $\alpha$-labeling then there exists
a cyclic $\G$-system of order $2en+1$ for any positive integer $n$.
Also, labeled graphs are models for a broad
range of applications, see for instance \cite{BG, BG1}.
For a very rich survey on graceful labelings we refer to \cite{G}.
Here we recall the basic definition. A \emph{graceful
labeling} of a graph $\G$ of size $e$ is an injective function $f: V(\G)\rightarrow
\{0,1,2,\ldots,e\}$ such that
$$\{|f(x)-f(y)| \ |\ [x,y]\in E(\G)\}=\{1,2,\ldots,e\}.$$
In the case where $\G$ is bipartite and $f$ has the additional property that its
maximum value on one of the two bipartite sets does not reach its minimum value on the other
one, one says that $f$ is an $\alpha$-\emph{labeling}.

Many variations of graceful labelings have been considered over the years. In particular Gnana
Jothi \cite{GJ} defines an \emph{odd graceful labeling} of a graph $\G$ of size $e$
as an injective function $f: V(\G)\rightarrow
\{0,1,2,\ldots,2e-1\}$ such that
$$\{|f(x)-f(y)| \ |\ [x,y]\in E(\G)\}=\{1,3,5,\ldots,2e-1\}.$$
In a recent paper, see \cite{APArs}, the second author have introduced the following new definition
which is, at the same time, a generalization of the concepts of a graceful labeling
(when $d=1$) and of an odd graceful labeling (when $d=e$).
\begin{defi}
Let $\G$ be a graph of size
$e=d\cdot m$. A $d$-\emph{divisible graceful labeling} of $\G$
is an injective function $f:V(\G) \rightarrow \{0,1,2,\ldots, d(m+1)-1\}$
such that
\begin{align*}
\{|f(x)-f(y)|\ |\ [x,y]\in E(\G)\} &=\{1,2,3,\ldots,d(m+1)-1\}\\
&\quad  \setminus \{m+1,2(m+1),\ldots,(d-1)(m+1)\}.
\end{align*}
Namely the set $\{|f(x)-f(y)|\ |\ [x,y]\in E(\G)\} $
can be divided into $d$ parts $P^0,P^1,\ldots,P^{d-1}$
where $P^i:=\{(m+1)i+1,(m+1)i+2,\ldots,(m+1)i+m\}$ for any
$i=0,1,\ldots, d-1$.
\end{defi}
\noindent
By saying that $d$ is \emph{admissible} we will mean that it is a divisor of $e$ and
so it makes sense to investigate the existence of a $d$-divisible graceful labeling of $\G$.
Also $\alpha$-labelings can be generalized in a similar way.
\begin{defi}
A $d$-\emph{divisible} $\alpha$-\emph{labeling} of a bipartite graph $\G$ is a
$d$-divisible graceful labeling of $\G$  having the property that
its maximum value on one of the two bipartite sets does not reach its minimum value on the other one.
\end{defi}

\noindent
Results on the existence of $d$-divisible
$\alpha$-labelings can be found in \cite{AP, APArs}.
In particular, using the notion
 of a $(v,d,\G,1)$-difference family introduced in \cite{BP},
 in \cite{APArs} it is proved that the existence of a $d$-divisible
($\alpha$-)labeling of a graph $\Gamma$ implies the existence of cyclic
$\Gamma$-decompositions. In fact:

\begin{thm}
If there exists a $d$-divisible graceful labeling of a graph $\G$ of size $e$ then there exists
a cyclic $\G$-decomposition of $K_{\left(\frac{e}{d}+1\right)\times 2d}$.
\end{thm}

\begin{thm}\label{decomp}
If there exists a $d$-divisible $\alpha$-labeling of a graph $\G$ of size $e$ then there exists
a cyclic $\G$-decomposition of $K_{(\frac{e}{d}+1)\times 2dn}$ for any positive integer $n$.
\end{thm}
\noindent
In this paper we deal with the existence of $d$-divisible $\alpha$-labelings of caterpillars,
hairy cycles and cycles. Often, in the literature, the problem of determining the existence of a labeling
of a given graph is solved by an ``ad hoc'' construction directly related to the case considered.
Here, in Section 2, we will construct a $d$-divisible $\alpha$-labeling of a caterpillar for any admissible value
of $d$. Also, we define a new kind of labeling for a graph $\G$, we show
that if $\G$ is a caterpillar such a labeling always exists and we call it \emph{standard $\alpha_S$-labeling} of $\G$.
 In Section 3, starting from a standard $\alpha_S$-labeling of a
suitable caterpillar
 we will construct a $d$-divisible $\alpha$-labeling of a hairy cycle or a cycle.
In detail, we will show that bipartite hairy cycles admit an odd $\alpha$-labeling
and, when all the vertices of the cycle
have the same degree, they admit a $d$-divisible $\alpha$-labeling
for any admissible value of $d$. About cycles we will show that,
for any positive integer $k$, $C_{4k}$ admits a
$d$-divisible $\alpha$-labeling for any admissible value of $d$.\\
We point out that in the cases of
caterpillars and cycles, if we consider the $d$-divisible $\alpha$-labelings here presented with $d=1$
we can obtain again the $\alpha$-labelings presented by Rosa in \cite{R}.\\
The existence of these $d$-divisible $\alpha$-labelings allows us to obtain
new infinite classes of decompositions of the complete multipartite graph.
 Here we obtain
a cyclic $\G$-decomposition of $K_{m\times n}$ for infinite pairs $(m,n)$
and for any caterpillar or hairy cycle $\G$.\\
The problem of  decomposing $K_{m\times n}$ into cycles has been taken in part;
for instance when the cycle length is a prime \cite{MP2006} , twice a prime \cite{S2008}
or ``small'' \cite{CB2000}.
In \cite{C1998, BCS2009, BCS2010} the authors gave necessary and sufficient condition for
the existence of cycle decompositions of  $K_{3\times n}$, $K_{4\times n}$
and $K_{5\times n}$, respectively. Here we obtain cyclic $C_{4k}$-decompositions
of  $K_{m\times n}$ for any positive integer $k$ and for infinite pairs $(m,n)$.

\section{Caterpillars}
\begin{figure}[H]
\begin{center}
\psfrag{x11}{$x_1^1$}
\psfrag{x1}{$x_1$}
\psfrag{x1n1}{$x_1^{n_1}$}
\psfrag{y11}{$y_1^1$}
\psfrag{y1}{$y_1$}
\psfrag{y1m1}{$y_1^{m_1}$}
\psfrag{x21}{$x_2^1$}
\psfrag{x2}{$x_2$}
\psfrag{x2n2}{$x_2^{n_2}$}
\psfrag{y21}{$y_2^1$}
\psfrag{y2}{$y_2$}
\psfrag{y2m2}{$y_2^{m_2}$}
\psfrag{x31}{$x_3^1$}
\psfrag{x3}{$x_3$}
\psfrag{x3n3}{$x_3^{n_3}$}
\psfrag{yt1}{$y_t^1$}
\psfrag{xt}{$x_t$}
\psfrag{yt}{$y_t$}
\psfrag{ytmt}{$y_t^{m_t}$}
\psfrag{1}{\tiny{$1$}}
\psfrag{n1}{\tiny{$n_1$}}
\psfrag{n1+1}{\tiny{$n_1+1$}}
\psfrag{n1+2}{\tiny{$n_1+2$}}
\psfrag{n1+m1+1}{\tiny{$n_1+m_1+1$}}
\psfrag{k}{\tiny{$k$}}
\includegraphics[width=0.8\textwidth]{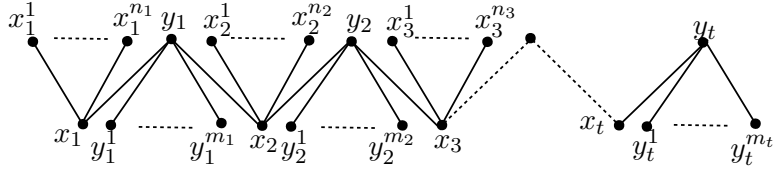}
\caption{The caterpillar $C[n_1,m_1\ldots,n_t,m_t]$}\label{Sunxy}
\end{center}
\end{figure}

\begin{defi}
A \emph{caterpillar} $\G$ is a tree, namely a graph without cycles,
with the property that the removal of its vertices of
degree one leaves a path $P$, called \emph{the path associated to} $\G$.
\end{defi}
\noindent
By the definition, any path is a caterpillar too.
In \cite{R}, Rosa proved that any caterpillar
admits an $\alpha$-labeling. In this section we will show, more in general,
that any caterpillar admits a $d$-divisible $\alpha$-labeling for any
admissible value of $d$. \\
First of all we introduce the following notation for the vertices of a caterpillar to emphasize that
we will see it as a bipartite graph.\\
Let $\G$ be a caterpillar and let $P$ be the path associated to $\G$.
If the length of $P$ is even,
set $P=[x_1,y_1,x_2,y_2,\ldots,x_t,y_t]$.
Let $n_i$ be the number of pendant edges through the vertex $x_i$
and $m_i$ be the number of pendant edges through the vertex $y_i$ for $i=1,\ldots,t$.
Also we call $x_i^1,x_i^2,\ldots,x_i^{n_i}$ the vertices of the edges pendant through
the vertex $x_i$ and $y_i^1,y_i^2,\ldots,y_i^{m_i}$ the vertices of the edges pendant through
the vertex $y_i$ for any $i=1,\ldots,t$; see Figure \ref{Sunxy}.
Hence the two bipartite sets are
\begin{eqnarray}
A=\{x_1,y_1^1,y_1^2,\ldots,y_1^{m_1},x_2,y_2^1,y_2^2,\ldots,y_2^{m_2},\ldots,x_t,y_t^1,y_t^2,\ldots,y_t^{m_t}\}\label{A} \\
B=\{x_1^1,x_1^2,\ldots,x_1^{n_1},y_1,x_2^1,x_2^2,\ldots,x_2^{n_2},y_2,\ldots,x_t^1,x_t^2,\ldots,x_t^{n_t},y_t\}\label{B}
\end{eqnarray}
and we will denote such a caterpillar by $C[n_1,m_1,\ldots,n_t,m_t]$.
\noindent
Notice that the elements of $E(\G)$ can be taken in a natural order from
$[x_1,x_1^1]$ up to $[y_t,y_t^{m_t}]$, that is we can write
\begin{eqnarray}
E(\G)=\{[x_1,x_1^1],\ldots,[x_1,x_1^{n_1}],[x_1,y_1],\ldots,[y_1,y_1^{m_1}],\cr
[y_1,x_2],\ldots,[x_t,y_t],[y_t,y_t^1],\ldots,[y_t,y_t^{m_t}]\}.\label{E}
\end{eqnarray}
Obviously if $n_1=0$,  (\ref{E}) becomes
\begin{eqnarray}\nonumber
E(\G)=\{[x_1,y_1],\ldots,[y_1,y_1^{m_1}],
[y_1,x_2],\ldots,[x_t,y_t],[y_t,y_t^1],\ldots,[y_t,y_t^{m_t}]\}.
\end{eqnarray}
Analogously, if $m_t=0$, we have
\begin{eqnarray}\nonumber
E(\G)=\{[x_1,x_1^1],\ldots,[x_1,x_1^{n_1}],[x_1,y_1],\ldots,[y_1,y_1^{m_1}],
[y_1,x_2],\ldots,[x_t,y_t]\}.
\end{eqnarray}
In what follows, by \emph{consecutive}, we will mean two vertices of $A$ or of $B$ or two edges of $E(\G)$ consecutive in the order
assumed in (\ref{A}), (\ref{B}) or (\ref{E}), respectively.\\
If the length of $P$ is odd, we can use analogous notations.\\



\noindent Now we introduce a new concept, which will be very useful in the following, that generalizes
that of a $d$-divisible $\alpha$-labeling.
\begin{defi}
Let $\G$ be a bipartite graph of size $e$ with parts $A$ and
$B$. Let $S$ be a set of $e$ positive integers.
We call \emph{$\alpha_S$-labeling of $\G$} an injective function $f:V(\G)\rightarrow \{0,1,2,\ldots,\max_{s\in S}s\}$
such that
$$\{|f(x)-f(y)|\ |\ [x,y]\in E(\G)\}=S\quad \textrm{and} \quad \mbox{$\max$}_A f<\mbox{$\min$}_B f.$$
\end{defi}

\begin{thm}\label{thm:alfaS}
Any caterpillar $\G$ of size $e$ admits an $\alpha_S$-labeling for any set $S$ of $e$ positive integers.
\end{thm}
\noindent
Proof.
Let $\G$ be a caterpillar with parts $A$ and $B$ defined as in (\ref{A}) and (\ref{B}), respectively.
We label the edges of $\G$ taken in the same order of (\ref{E}),
with the elements of $S$ in ascending order. Let $f:V(\G)\rightarrow \{0,1,2,\ldots,\max_{s\in S}s\}$
be the function defined so that,
given
$x\in A$ and $y \in B$, the label of the edge $[x,y]$ is $f(y)-f(x)$ and with $f(y_t^{m_t})=0$
$(f(x_t)=0$ if $m_t=0)$.
Note that if we consider the elements of $A$ in the same order of (\ref{A})
their labels are in descending order and if we consider the elements of $B$ in the same order of (\ref{B})
their labels are in ascending order.
So, it is easy to see that $f$ is an injective function and that $\max_A f< \min_B f$.
Hence $f$ is an $\alpha_S$-labeling of $\G$.
\hfill $\Box$

\begin{defi}\label{OnS}
We will call the function $f$ constructed in the proof of Theorem
\ref{thm:alfaS} \emph{the standard $\alpha_S$-labeling of $\G$}.
\end{defi}

\noindent
Thanks to Theorem \ref{thm:alfaS} we can completely solve the problem of the existence
of $d$-divisible $\alpha$-labelings of caterpillars.
\begin{cor}\label{caterpillar}
Any caterpillar admits a $d$-divisible $\alpha$-labeling for any admissible value of $d$.
\end{cor}
\noindent
Proof. Let $\G$ be a caterpillar with $e=d\cdot m$ edges.
Let $f$ be the standard $\alpha_S$-labeling of $\G$ where
$S=\{1,2,\ldots,d(m+1)-1\}\setminus\{m+1,2(m+2),\ldots,(d-1)(m+1)\}$.
It is easy to see that $f$
is a $d$-divisible $\alpha$-labeling of $\G$.
\hfill $\Box$\\
\\

\noindent
As an immediate consequence of Theorem \ref{decomp} and Corollary \ref{caterpillar} we have:
\begin{thm}
Let $\G$ be a caterpillar with $e$ edges. There exists a cyclic $\G$-decomposition of $K_{(\frac{e}{d}+1)\times 2dn}$
for any divisor $d$ of $e$ and any positive integer $n$.
\end{thm}

\begin{rem}
We point out that if $d=1$ the $d$-divisible $\alpha$-labeling of the proof of Corollary \ref{caterpillar}
is nothing but the $\alpha$-labeling obtained by Rosa. Also,
if the caterpillar is indeed a path,
we find again the $d$-divisible $\alpha$-labeling constructed in
\cite{APArs}.
\end{rem}
\begin{ex}
We consider a caterpillar $\G$ with $12$ edges.
In Figure \ref{ExCater} we show all the possible $d$-graceful $\alpha$-labelings of $\G$,
other than the classical one, obtained following
the proof of Corollary \ref{caterpillar}.
\begin{figure}[H]
\begin{center}
\psfrag{a}{$\alpha$}
\psfrag{b}{$\beta$}
\psfrag{0}{\small{$0$}}
\psfrag{1}{\small{$1$}}
\psfrag{2}{\small{$2$}}
\psfrag{3}{\small{$3$}}
\psfrag{4}{\small{$4$}}
\psfrag{5}{\small{$5$}}
\psfrag{6}{\small{$6$}}
\psfrag{7}{\small{$7$}}
\psfrag{8}{\small{$8$}}
\psfrag{9}{\small{$9$}}
\psfrag{10}{\small{$10$}}
\psfrag{11}{\small{$11$}}
\psfrag{12}{\small{$12$}}
\psfrag{13}{\small{$13$}}
\psfrag{14}{\small{$14$}}
\psfrag{15}{\small{$15$}}
\psfrag{16}{\small{$16$}}
\psfrag{17}{\small{$17$}}
\psfrag{18}{\small{$18$}}
\psfrag{19}{\small{$19$}}
\psfrag{20}{\small{$20$}}
\psfrag{21}{\small{$21$}}
\psfrag{22}{\small{$22$}}
\psfrag{23}{\small{$23$}}
\psfrag{s1}{\tiny{$1$}}
\psfrag{s2}{\tiny{$2$}}
\psfrag{s3}{\tiny{$3$}}
\psfrag{s4}{\tiny{$4$}}
\psfrag{s5}{\tiny{$5$}}
\psfrag{s6}{\tiny{$6$}}
\psfrag{s7}{\tiny{$7$}}
\psfrag{s8}{\tiny{$8$}}
\psfrag{s9}{\tiny{$9$}}
\psfrag{s10}{\tiny{$10$}}
\psfrag{s11}{\tiny{$11$}}
\psfrag{s12}{\tiny{$12$}}
\psfrag{s13}{\tiny{$13$}}
\psfrag{s14}{\tiny{$14$}}
\psfrag{s15}{\tiny{$15$}}
\psfrag{s16}{\tiny{$16$}}
\psfrag{s17}{\tiny{$17$}}
\psfrag{s18}{\tiny{$18$}}
\psfrag{s19}{\tiny{$19$}}
\psfrag{s20}{\tiny{$20$}}
\psfrag{s21}{\tiny{$21$}}
\psfrag{s22}{\tiny{$22$}}
\psfrag{s23}{\tiny{$23$}}
\subfigure[]{
$\begin{array}{c}
\includegraphics[width=0.15\textwidth]{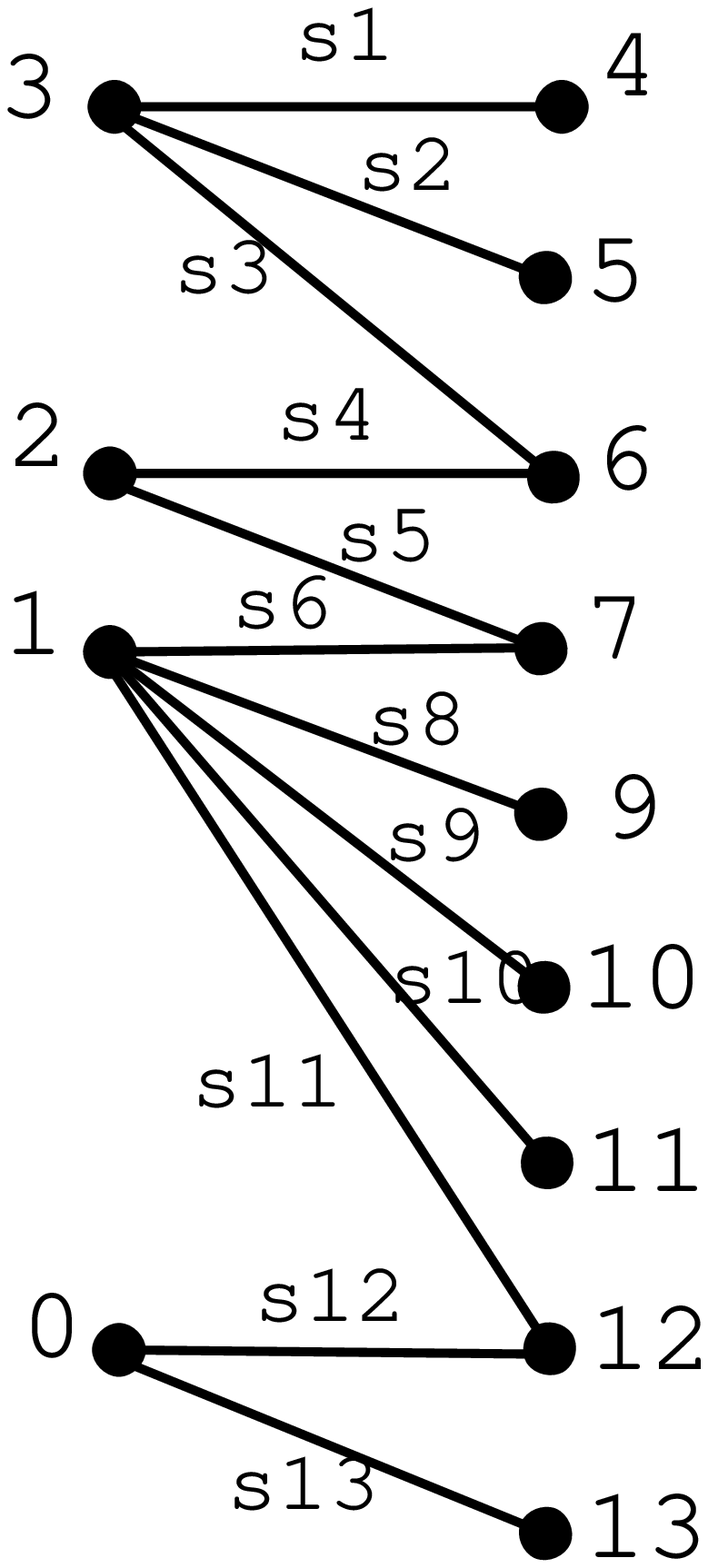}
\end{array}$}
\subfigure[]{
$\begin{array}{c}
\includegraphics[width=0.15\textwidth]{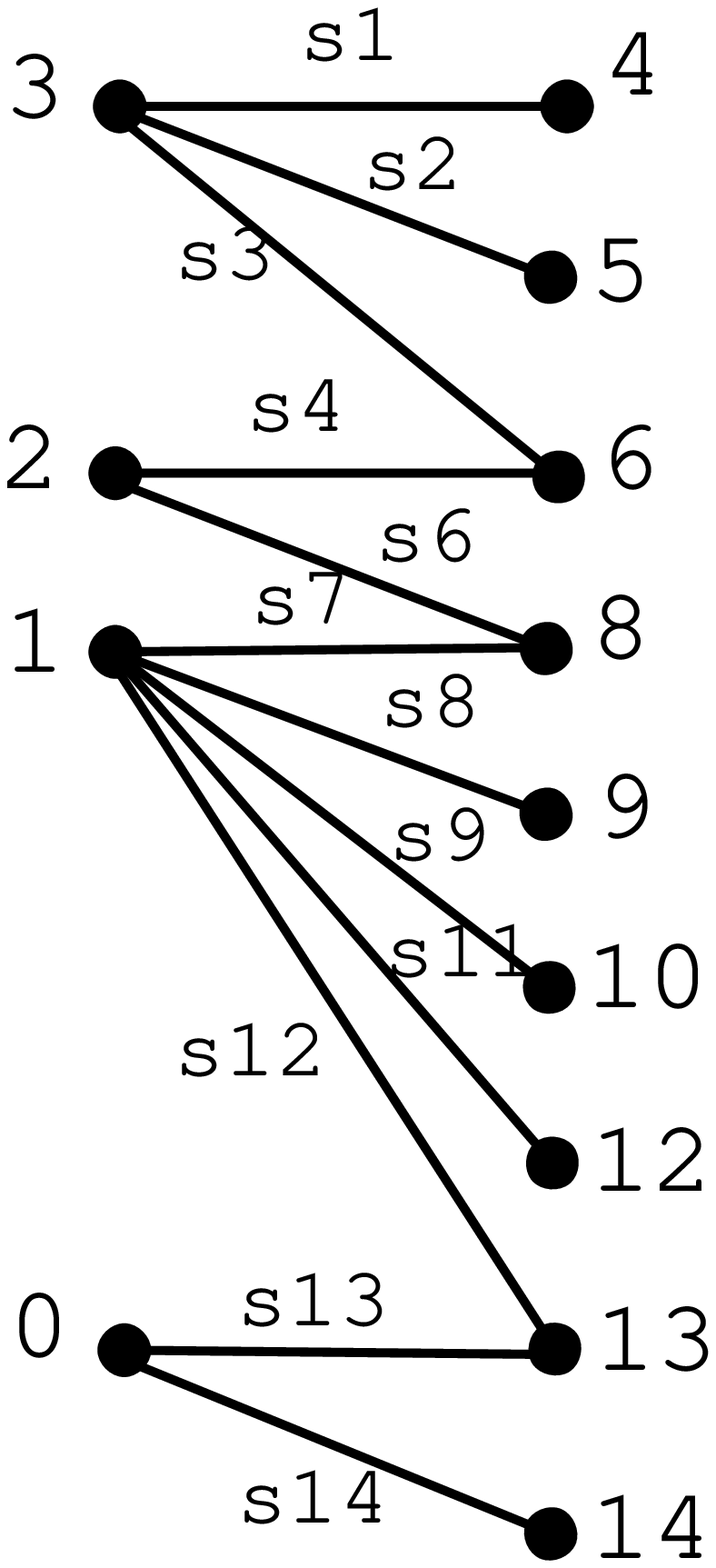}
\end{array}$}
\subfigure[]{
$\begin{array}{c}
\includegraphics[width=0.15\textwidth]{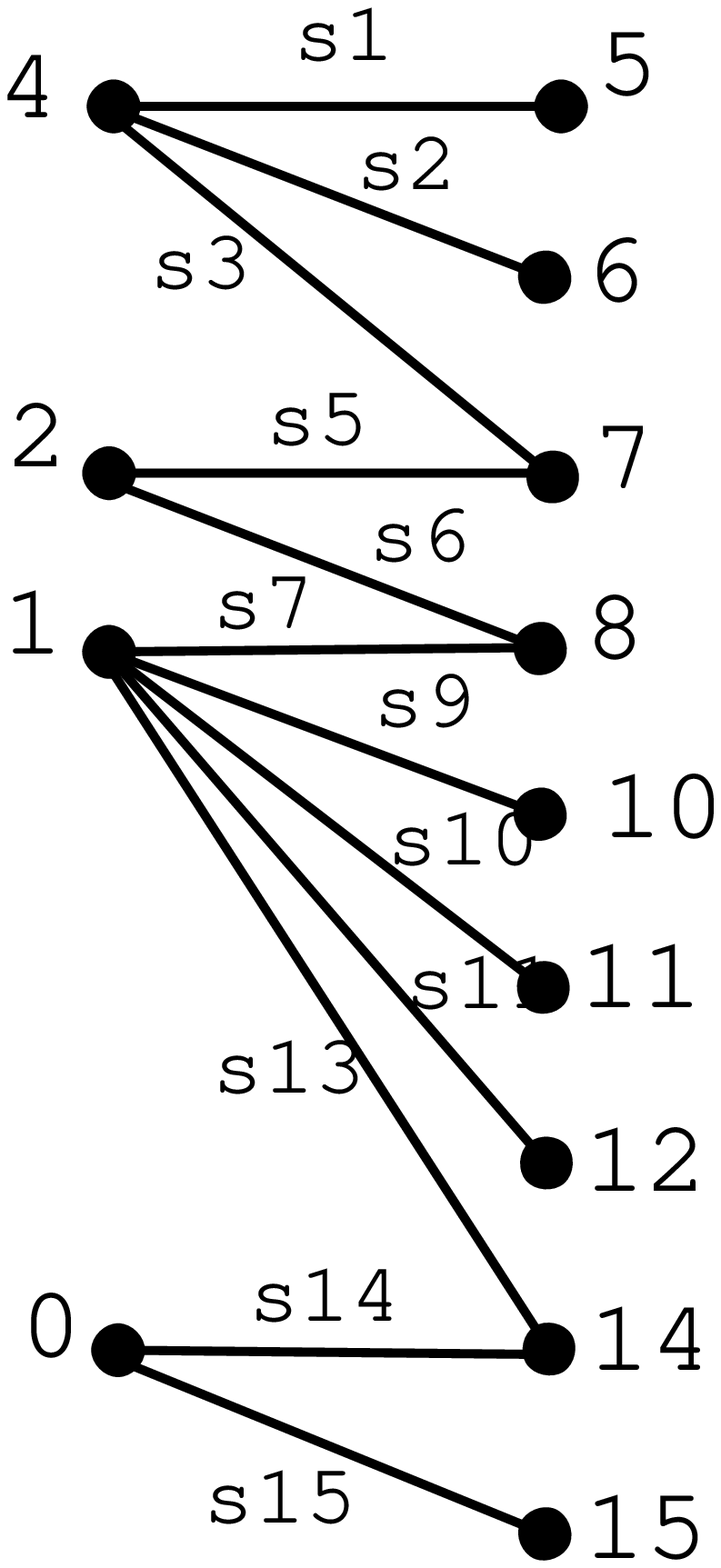}
\end{array}$}
\subfigure[]{
$\begin{array}{c}
\includegraphics[width=0.15\textwidth]{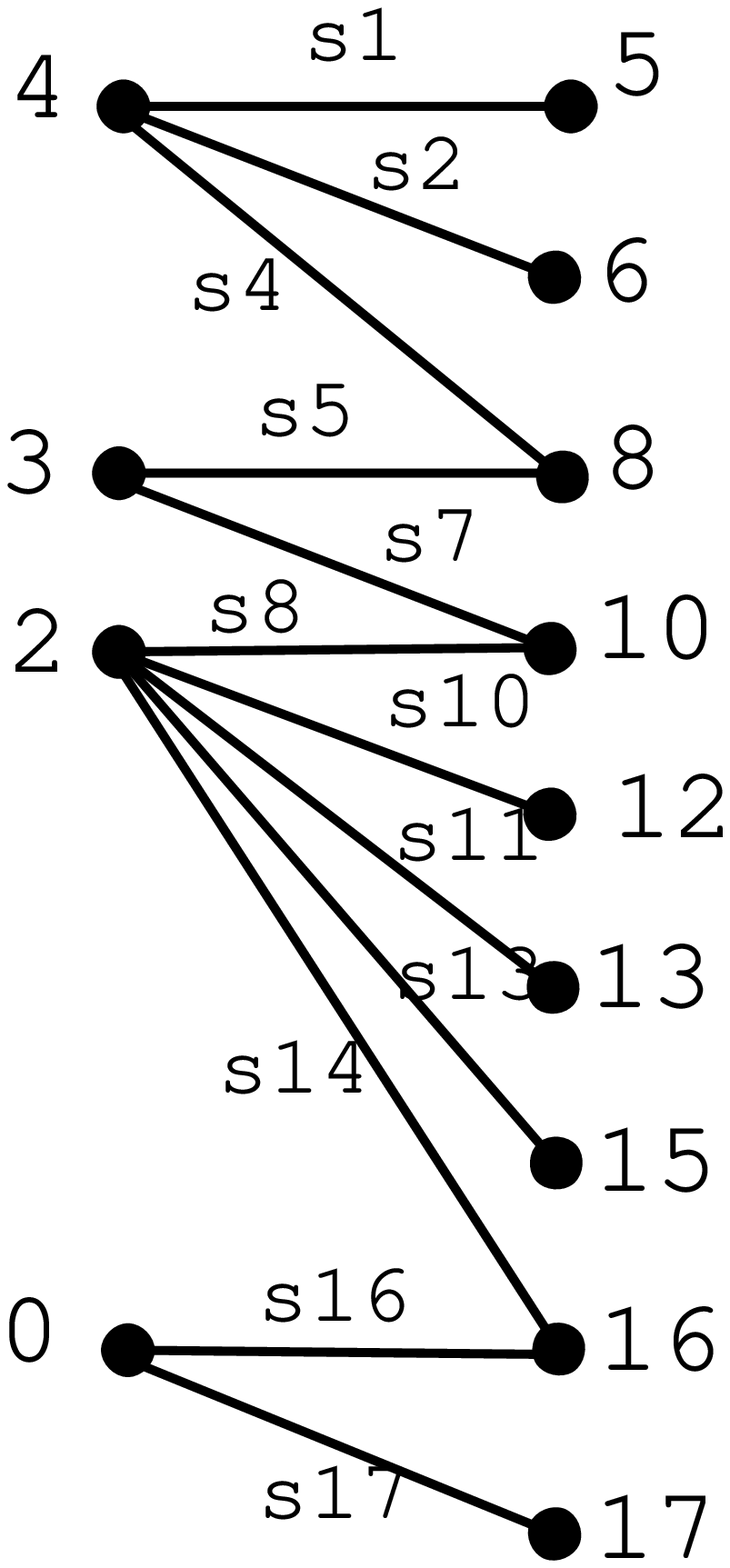}
\end{array}$}
\subfigure[]{
$\begin{array}{c}
\includegraphics[width=0.15\textwidth]{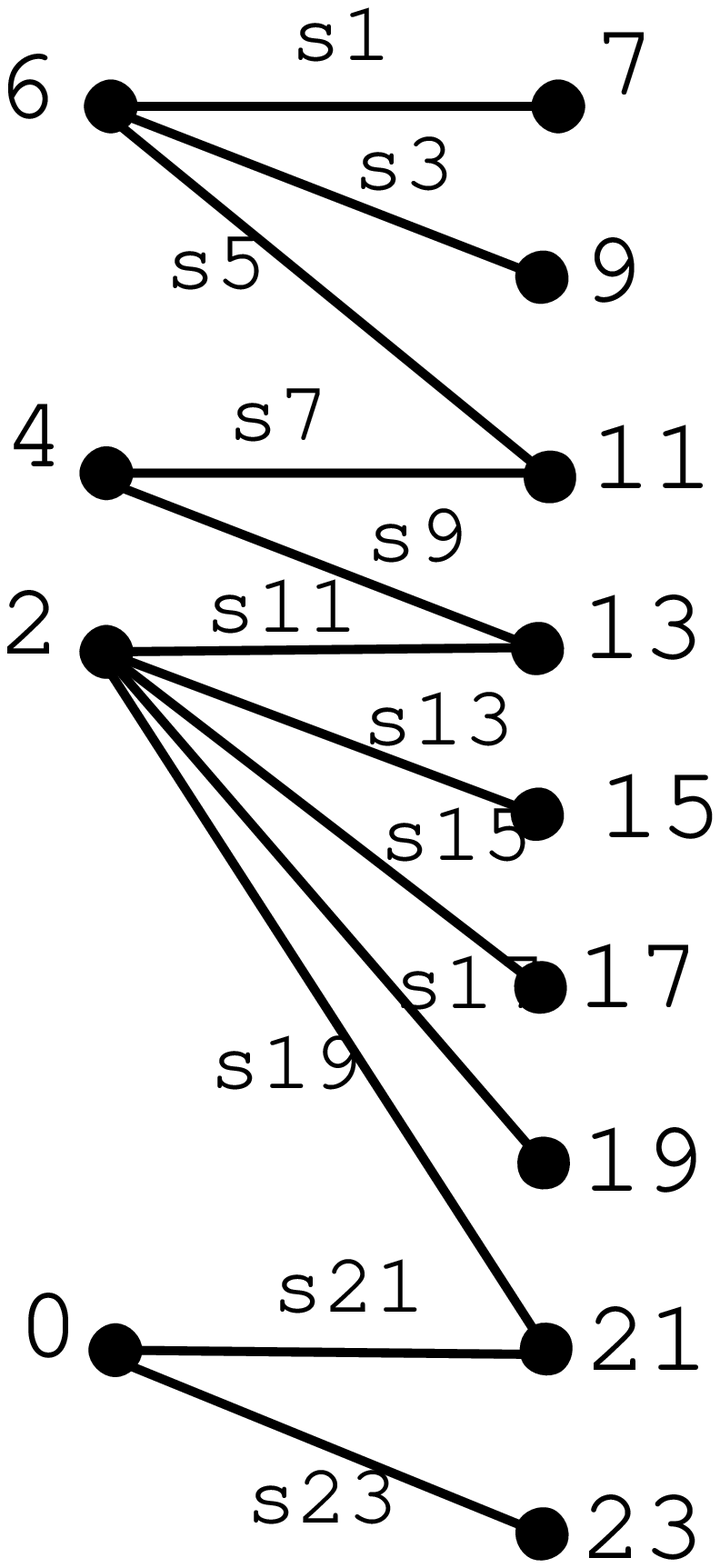}
\end{array}$}
\caption{$d$-divisible $\alpha$-labelings of a caterpillar}\label{ExCater}
\end{center}
\end{figure}
\end{ex}

\noindent
\begin{defi}
Let $f$ be an $\alpha_S$-labeling of a graph $\G$ and let $X\subseteq V(\G)$.
We will call \emph{missing vertex label} in $f(X)$
(mv-label, for short) each element of $\{n\in \mathbb{N}\ \mid\ \min_X f\leq n\leq \max_X f\}\setminus f(X)$.
\end{defi}

\noindent
We point out that following the construction of the
$\alpha_S$-labeling provided in Theorem \ref{thm:alfaS} each missing integer in $S$ causes a
corresponding mv-label in $f(A)\cup f(B)$. For instance, look at Figure \ref{ExCater}(d).
Now $S=\{1,2,\ldots,17\}\setminus\{3,6,9,12,15\}$ and the missing integers
$\{3,6,9,12,15\}$ cause the mv-labels $\{7,9,11,14,1\}$,
respectively.


\section{Hairy cycles  and cycles from caterpillars}
\begin{defi}
A unicyclic graph $\Gamma$, other than a cycle, is called a \emph{hairy cycle} if
the deletion of any edge in the cycle of $\G$ results in a caterpillar.
\end{defi}
So, both cycles and hairy cycles can be always seen as suitable caterpillars $\G$ with one extra edge:
the one connecting the  ending vertices of the path associated to $\G$.

Following this line and keeping in mind the construction of Rosa (see \cite{R}), in \cite{Ba} Barrientos proves that
all hairy cycles are graceful and that when the graph is bipartite, namely when the cycle has even length, the
labeling is an $\alpha$-labeling.

\begin{figure}[H]
\begin{center}
\psfrag{a}{$\alpha$}
\psfrag{b}{$\beta$}
\psfrag{yt}{$y_t$}
\psfrag{x1}{$x_1$}
\subfigure[]{
$\begin{array}{c}
\includegraphics[width=0.15\textwidth]{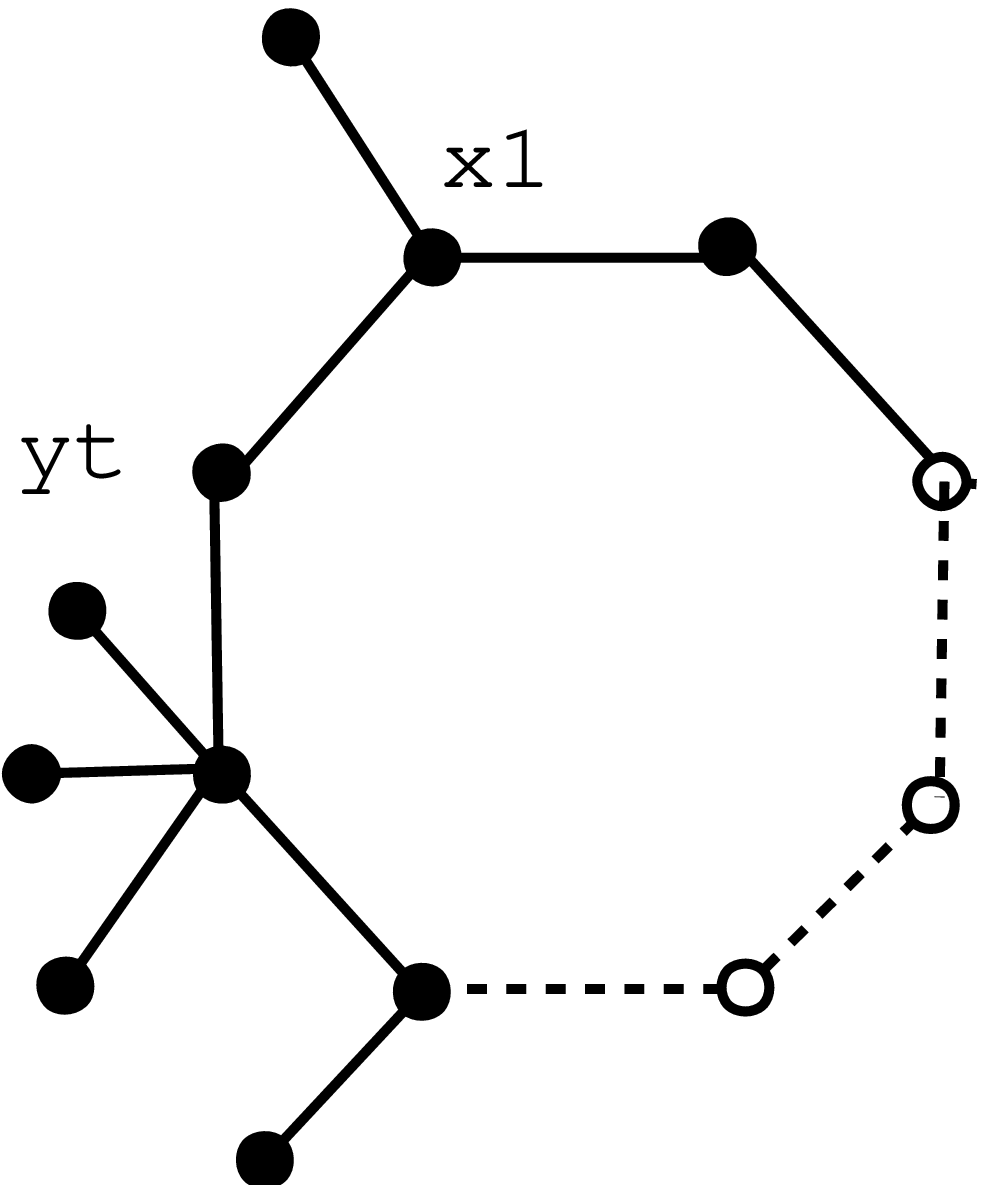}\\
\\
\includegraphics[width=0.25\textwidth]{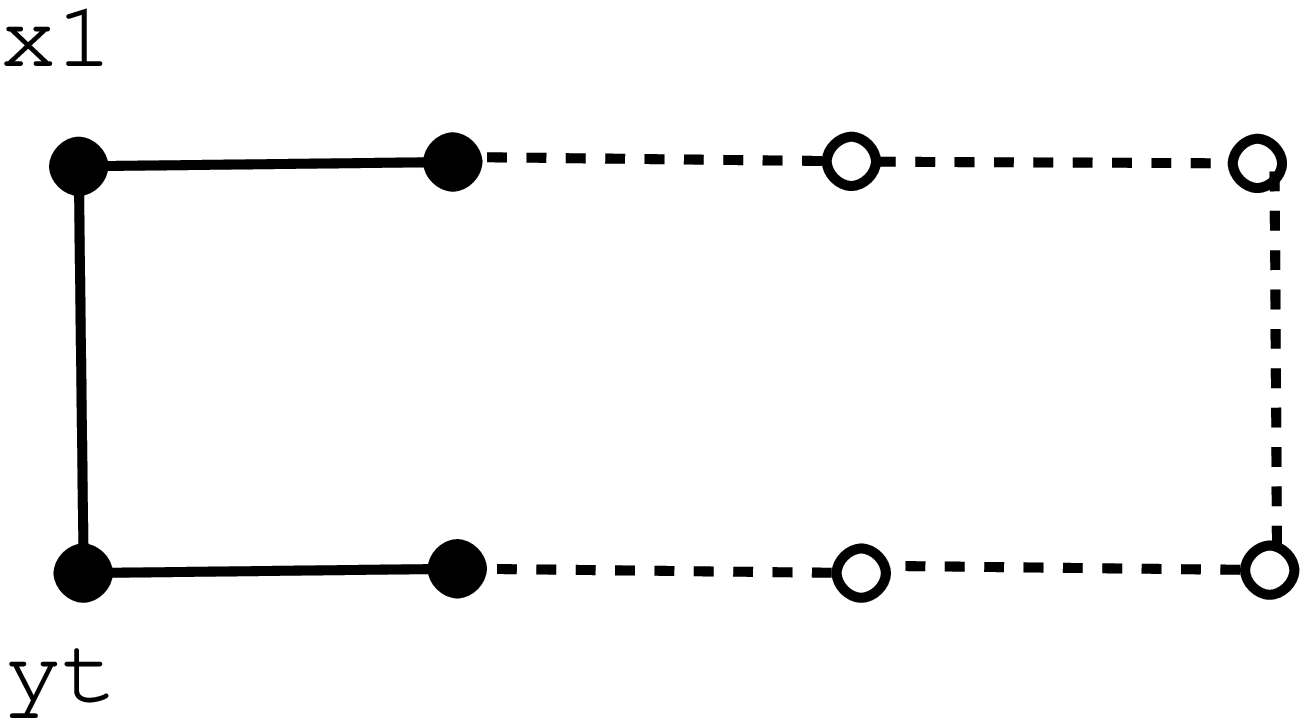}
\end{array}$}
\quad\quad\quad\quad
\subfigure[]{
$\begin{array}{c}
\vspace{0,7cm}
\includegraphics[width=0.4\textwidth]{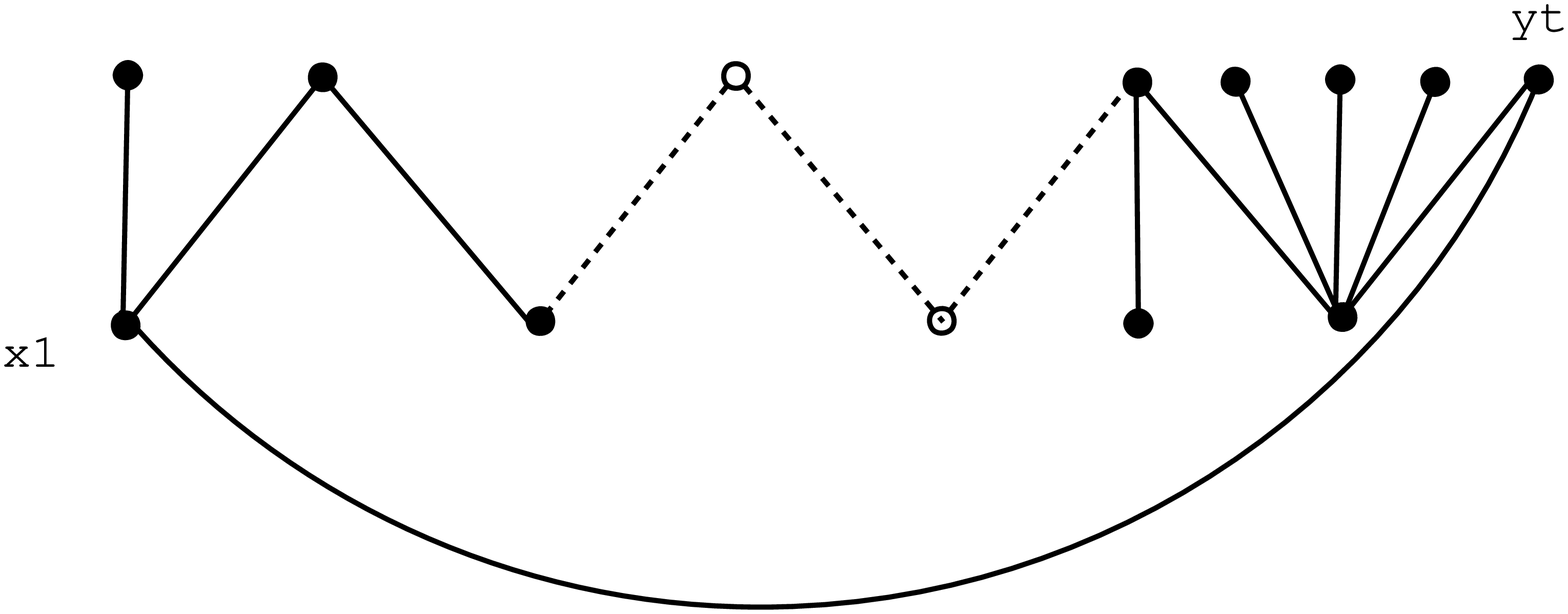}\\
\includegraphics[width=0.4\textwidth]{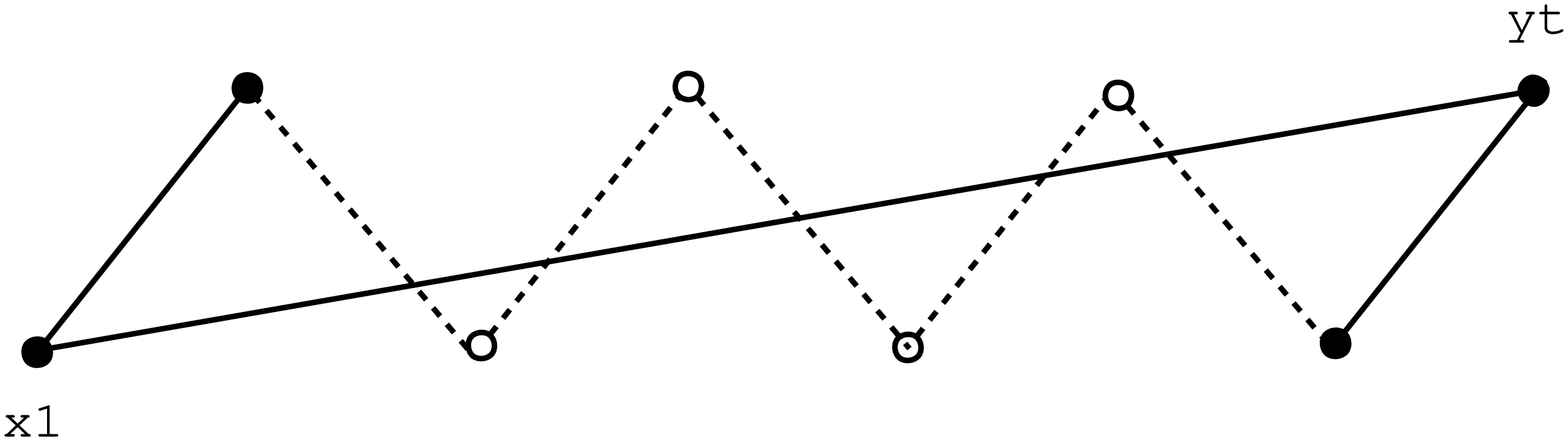}
\end{array}$}
\caption{An hairy cycle and a cycle}\label{Sun}
\end{center}
\end{figure}

Following the same line and using the results of Section 2, in the sequel we will give a method to construct
$d$-divisible $\alpha$-labelings of hairy cycles and cycles.
Thus, from now on we will draw hairy cycles and cycles  as in Figure \ref{Sun}(b) but instead
in the classical way of Figure \ref{Sun}(a) and by $HC(n_1,m_1,\dots,n_t,m_t)$, with $t\geq2$, we will denote the hairy cycle obtainable from
the caterpillar $C[n_1,m_1,\dots,n_t,m_t]$ by adding the edge $[x_1,y_t]$.
 If the cycle $(x_1,y_1,\dots,x_{t-1},y_{t-1},x_t)$, with $t\geq2$, had odd length, we would obviously refer
for the constructed hairy cycle with
$HC(n_1,m_1,\dots,n_{t-1},m_{t-1},n_t)$.

\begin{rem}\label{rem:c}
Consider  the graphs of Figure \ref{Sun}(b) and let $e$ be their size. Removing  the edge $[x_1,y_t]$
we obtain, in both cases, a caterpillar $\Gamma$ of size $e-1$. Let $S=\{a_1,a_2,\dots,a_{e}\}$ be a set of $e$ integers.
Set $S'=\{a_1,a_2,\dots,a_{e}\}\setminus\{a_c\}$, $1\leq c \leq e$, and consider an
$\alpha_{S'}$-labeling of $\Gamma$, say $f$. If it results in $f(y_t)-f(x_1)=a_c$,
 in a natural way we can extend $f$ to an $\alpha_S$-labeling of the original graph,
 labeling the edge $[x_1,y_t]$ by $a_c$.
\end{rem}
\noindent
Let the elements of $S$ be taken in ascending order. Let $S'$ be defined as in previous remark and $f$ be the standard $\alpha_{S'}$-labeling of $\G$.
Obviously, the positions of the mv-labels affect the value of $f(y_t)-f(x_1)$.
In particular, notice that if $a_{c-1}$ and $a_{c+1}$ label two consecutive edges through a vertex of $B$ $($or $A)$,
the deletion of $a_c$ causes an
extra mv-label in $f(A)$ $($or $f(B))$.\\
So, given an $\alpha_S$-labeling $f$ of a caterpillar $\G$ we can try to lightly modify $f$, if it is the case,
changing the position of some mv-labels in order
to obtain another function $g$ which results in still an $\alpha_S$-labeling of $\G$
with $g(y_t)-g(x_1)\neq f(y_t)-f(x_1)$. Here we are giving some few ways to obtain functions as $g$.\\
In the sequel $\G$ is a caterpillar
(see Figure \ref{Sunxy}), $S$ is an arbitrary set of positive integers of size $e=|E(\G)|$,
$f$ is the standard $\alpha_S$-labeling of $\G$. If $x$ and $y$ are consecutive vertices of $A$ or of $B$,
$f(x)+h=f(y)$ means that there are $h-1$ mv-labels between $f(x)$ and $f(y)$.\\
\\
\noindent
$\bullet$~$\mathbf{[O_1]}$~~gives $g(y_t)-g(x_1)=f(y_t)-f(x_1)-1$.\\
\emph{
If there exists $s\in \{1,2,\ldots,t-1\}$  such that the following conditions are true:\\
(1) $f(y_s)+1=f(x_{s+1}^1)$, \ \ \ (2) $f(y_s)=f(x_s^{n_s})+h$, with $h\geq 2$,\\
then, the function $g$ defined as follows:\\
$g(x_{s+1}^1)=f(x_{s+1}^1)-1$\\
$g(z)=f(z) \ \ \ \forall z\in \{y_t^{m_t},\dots,x_{s+1}\}\bigcup\{x_{s+1}^2,\dots,y_t\}=H_1$\\
$g(z)=f(z)+1 \ \ \ \forall z\in \{y_s^{m_s},\dots,x_{1}\}\bigcup\{x_{1}^1,\dots,y_s\}=H_2$\\
results in again an $\alpha_S$-labeling of $\Gamma$.}

In fact,\\
1.~$S$ results in the same because the labels of the edges from $[x_1,x_1^1]$ up to $[y_s^{m_s},y_s]$ as well as those
from $[x_{s+1},x_{s+1}^2]$ up to $[y_t^{m_t},y_t]$ are the same
while those of $[x_{s+1},y_s]$ and $[x_{s+1},x^1_{s+1}]$ commute.\\
2.~$g$ is injective because $g(z_1)\neq g(z_2)$ for $z_1,z_2\in H_1\bigcup H_2 \ \mbox{with} \ z_1\neq z_2$.
If, by contradiction, $g(y_s)=g(x^1_{s+1})$ then $f(y_s)+1= f(x^1_{s+1})-1$ that is $f(x^1_{s+1})=f(y_s)+2$,
excluded by (1). If, again by contradiction, $g(x^1_{s+1})=g(x_s^{n_s})$ then  $f(x^1_{s+1})-1=f(x_s^{n_s})+1$
that is $f(y_s)=f(x_s^{n_s})+1$, excluded by (2).\\
3.~$\mbox{max}_A~g=\mbox{max}_A~f+1<\mbox{min}_B~f+1=\mbox{min}_B~g$, so $g$ results in an $\alpha_S$-labeling of $\G$.\\

\noindent
$\bullet$~$\mathbf{[O_2]}$~~gives $g(y_t)-g(x_1)= f(y_t)-f(x_1)-2$.\\
\emph{
If there exists $s\in \{1,2,\ldots,t-1\}$ such that the following conditions are true:\\
(1) $f(y_s)+2=f(x_{s+1}^1)$, \ \ \ (2) $f(y_s)=f(x_s^{n_s})+h$, with $h\geq 3$,\\
then, the function $g$ defined as follows:\\
$g(x_{s+1}^1)=f(x_{s+1}^1)-2$\\
$g(z)=f(z) \ \ \ \forall z\in \{y_t^{m_t},\dots,x_{s+1}\}\bigcup\{x_{s+1}^2,\dots,y_t\}=H_1$\\
$g(z)=f(z)+2 \ \ \ \forall z\in \{y_s^{m_s},\dots,x_{1}\}\bigcup\{x_{1}^1,\dots,y_s\}=H_2$\\
results in again an $\alpha_S$-labeling of $\Gamma$.}

In fact,\\
1.~$S$ results in the same because the labels of the edges from $[x_1,x_1^1]$ up to $[y_s^{m_s},y_s]$ as well as those
from $[x_{s+1},x_{s+1}^2]$ up to $[y_t^{m_t},y_t]$ are the same
while those of $[x_{s+1},y_s]$ and $[x_{s+1},x^1_{s+1}]$ commute.\\
2.~$g$ is injective because $g(z_1)\neq g(z_2)$ for $z_1,z_2\in H_1\bigcup H_2 \ \mbox{with} \ z_1\neq z_2$.
If, by contradiction, $g(y_s)=g(x^1_{s+1})$ then $f(y_s)+2= f(x^1_{s+1})-2$ that is $f(x^1_{s+1})=f(y_s)+4$,
excluded by (1). If, again by contradiction, $g(x^1_{s+1})=g(x_s^{n_s})$ then  $f(x^1_{s+1})-2=f(x_s^{n_s})+2$
that is $f(y_s)=f(x_s^{n_s})+2$, excluded by (2).\\
3.~$\mbox{max}_A~g=\mbox{max}_A~f+2<\mbox{min}_B~f+2=\mbox{min}_B~g$, so $g$ results in an $\alpha_S$-labeling of $\G$.\\

\noindent
$\bullet$~$\mathbf{[O_3]}$~~gives $g(y_t)-g(x_1)= f(y_t)-f(x_1)-2$.\\
\emph{
If there exists $s\in \{1,2,\ldots,t-1\}$ such that the following conditions are true:\\
(1)~$f(y_s)+2=f(x_{s+1}^1)+1=f(x_{s+1}^2)$, \ \ \ (2)~$f(y_s)=f(x_s^{n_s})+h$,
with $h\geq 3$,\\
then, the function $g$ defined as follows:\\
$g(x_{s+1}^2)=f(x_{s+1}^2)-2$, $f(x^1_{s+1})=g(x^1_{s+1})$\\
$g(z)=f(z) \ \ \ \forall z\in \{y_t^{m_t},\dots,x_{s+1}\}\bigcup\{x_{s+1}^3,\dots,y_t\}=H_1$  \\
$g(z)=f(z)+2 \ \ \ \forall z\in \{y_s^{m_s},\dots,x_{1}\}\bigcup\{x_{1}^1,\dots,y_s\}=H_2$\\
results in again an $\alpha_S$-labeling of $\Gamma$.}

In fact,\\
1.~$S$ results in the same because the labels of the edges from $[x_1,x_1^1]$ up to $[y_s^{m_s},y_s]$, those of the edges
from $[x_{s+1},x_{s+1}^3]$ up to $[y_t^{m_t},y_t]$ as well as the label of  $[x_{s+1},x_{s+1}^1]$ are the same
while those of $[x_{s+1},y_s]$ and $[x_{s+1},x^2_{s+1}]$ commute.\\
2.~$g$ is injective because $g(z_1)\neq g(z_2)$ for $z_1,z_2\in H_1\bigcup H_2 \  \mbox{with} \ z_1\neq z_2$.
If, by contradiction, $g(y_s)=g(x^1_{s+1})$ then $f(y_s)+2= f(x^1_{s+1})$, excluded by (1).
If, by contradiction, $g(y_s)=g(x^2_{s+1})$ then $f(y_s)+2= f(x^2_{s+1})-2=f(y_s)$, excluded.
If, again by contradiction, $g(x^2_{s+1})=g(x_s^{n_s})$ then  $f(x^2_{s+1})-2=f(x_s^{n_s})+2$
that is $f(y_s)=f(x_s^{n_s})+2$, excluded by (2).\\
3.~$\mbox{max}_A~g=\mbox{max}_A~f+2<\mbox{min}_B~f+2=\mbox{min}_B~g$, so $g$ results in an $\alpha$-labeling of $\G$.\\

\noindent
$\bullet$~$\mathbf{[O_4]}$~~gives $g(y_t)-g(x_1)= f(y_t)-f(x_1)-h$ where $h=f(y_t)-f(x_t^{n_t})$.\\
\emph{
If the following conditions are true:    \ \ (1) $n_t\neq 0$, \ \ \ (2) $m_t=0$,\\
then, the function $g$ defined as follows:\\
$g(z)=f(z)\ \  \forall z\in V(\G)\setminus \{y_t,x_{t}^{n_t}\}$,\
$g(y_t)=f(x_{t}^{n_t})$ and $g(x_{t}^{n_t})=f(y_t)$\\
results in again an $\alpha_S$-labeling of $\Gamma$.}

In fact,\\
1.~$S$ results in the same because the labels of the edges are the same with the exception of
 those of $[x_t,y_t]$ and $[x_t,x^{n_t}_{t}]$ that commute.\\
2.~$g$ is injective because  $g(z)=f(z)$ for each $z\in V(\G)$ with the exception of $y_t$ and $x_{t}^{n_t}$,
whose images are swapped.\\
3.~$\mbox{max}_A~g=\mbox{max}_A~f<\mbox{min}_B~f=\mbox{min}_B~g$, so $g$ results in an $\alpha_S$-labeling of $\G$.\\

\noindent$\bullet$~$\mathbf{ [O_5]}$~~gives $g(y_t)-g(x_1)= f(y_t)-f(x_1)-1$.\\
\emph{
If there exists $s\in \{1,2,\ldots,t-1\}$ such that the following conditions are true:\\
(1)~$f(y_s)=f(x_{s+1}^1)-2$, \ \ \ (2)~$f(y_s)=f(x_s^{n_s})+1$,\\
then, the function $g$ defined as follows:\\
$g(x_{s}^{n_s})=f(x_{s}^{n_s})+2$, \\
$g(z)=f(z) \ \ \ \forall z\in \{y_t^{m_t},\dots,y_s^1\}\bigcup\{y_s,\dots,y_t\}=H_1$  \\
$g(z)=f(z)+1 \ \ \ \forall z\in \{x_s,\dots,x_{1}\}\bigcup\{x_{1}^1,\dots,x_s^{n_s-1}\}=H_2$\\
results in again an $\alpha_S$-labeling of $\Gamma$.}

In fact,\\
1.~$S$ results in the same because the labels of the edges from $[x_1,x_1^1]$ up to $[x_s^{n_s-1},x_s]$ and
those of the edges from $[y_{s},y_{s}^1]$ up to $[y_t^{m_t},y_t]$ are the same
while those of $[x_{s},x_s^{n_s}]$ and $[x_{s},y_s]$ commute.\\
2.~$g$ is injective because $g(z_1)\neq g(z_2)$ for $z_1,z_2\in H_1\bigcup H_2
 \  \mbox{with} \ z_1\neq z_2$.
If, by contradiction, $g(x_s^{n_s-1})=g(x_s^{n_s})$ then ${f}(x_s^{n_s-1})+1={f}(x_s^{n_s})+2$,
so ${f}(x_s^{n_s-1})={f}(x_s^{n_s})+1=f(y_s)$ (by (2)),  excluded as $f$ is injective.
Also, if, by contradiction, $g(x_s^{n_s})=g(x_{s+1}^{1})$ then
${f}(x_s^{n_s})+2={f}(x_{s+1}^{1})$ so $f(y_s)+1=f(y_s)+2$, by (1)
and (2), which is obviously impossible.\\
3.~$\mbox{max}_A~g=\mbox{max}_A~f+1<\mbox{min}_B~f+1=\mbox{min}_B~g$, so $g$ results in an
$\alpha_S$-labeling of $\G$.\\

\noindent$\bullet$~$\mathbf{ [O_5]_4}$~~gives $g(y_t)-g(x_1)= f(y_t)-f(x_1)-1$.\\
\emph{
If there exists $s\in \{1,2,\ldots,t-1\}$ such that the following conditions are true:\\
(1)~$\exists j\in\{1,2,\ldots,m_s\}$ such that $f(y_s^j)=f(x_s)-4$, \ \ \ (2)~$f(y_s)=f(x_s^{n_s})+1$,
\ (3)~$f(x_{s+1}^i)=f(y_s)+i, \ i=1,2,3,4$, and
$f(x_{s+1}^5)=f(y_s)+6$,\\
then, the function $g$ defined as follows:\\
$g(x_{s}^{n_s})=f(y_{s})+5$,  \ \ \ $g(y_{s}^{j})=f(x_{s})$\\
$g(z)=f(z) \ \ \ \forall z\in (\{y_t^{m_t},\dots,y_s^1\}\setminus\{y_{s}^{j}\})\bigcup\{y_s,\dots,y_t\}=H_1$  \\
$g(z)=f(z)+1 \ \ \ \forall z\in
\{x_s,\dots,x_{1}\}\bigcup\{x_{1}^1,\dots,x_s^{n_s-1}\}=H_2$\\
results in again an $\alpha_S$-labeling of $\Gamma$.}

In fact,\\
1.~$S$ results in the same because the labels of the edges from $[x_1,x_1^1]$ up to $[x_s^{n_s-1},x_s]$ and
those of the edges from $[y_{s},y_{s}^1]$ up to $[y_t^{m_t},y_t]$, excepted $[y_{s},y_{s}^j]$, are the same
while those of $[x_{s},x_s^{n_s}]$, $[x_{s},y_s]$ and $[y_{s},y_{s}^j]$ cyclically permute.\\
2.~$g$ is injective because $g(z_1)\neq g(z_2)$ for $z_1,z_2\in H_1\bigcup H_2
 \  \mbox{with} \ z_1\neq z_2$. If, by contradiction, there was an $i$ such that
 $g(x_{s}^{n_s})=g(x_{s+1}^i)$, then $f(y_{s})+5=f(x_{s+1}^i)$, excluded by (3).
If, by contradiction, $g(x_s^{n_s-1})=g(x_s^{n_s})$ then ${f}(x_s^{n_s-1})+1={f}(y_s)+5$,
so ${f}(x_s^{n_s-1})={f}(y_s)+4=f(x_{s+1}^4)$ (by (3))  excluded as $f$ is injective.\\
3.~$\mbox{max}_A~g=\mbox{max}_A~f+1<\mbox{min}_B~f+1=\mbox{min}_B~g$, so $g$ results in an
$\alpha_S$-labeling of $\G$.\\

\begin{rem}
Obviously, we can apply more than one of the previous operations
to the same standard $\alpha_S$-labeling, as long as they operate
on disjoint set of edges.
\end{rem}

\subsection{Hairy cycles}
In this subsection we focus our attention on
hairy cycles.
\noindent
In \cite{Ba} Barrientos gives a labeling for any hairy cycle and when the graph is
bipartite such a labeling is an $\alpha$-labeling.\\
Here we show that  any bipartite hairy cycle $H$ admits an odd $\alpha$-labeling
(namely an $e$-divisible $\alpha$-labeling, $e$ being the size of $H$).
Then, for any admissible value of $d$, we will prove the existence of $d$-divisible $\alpha$-labelings
of an infinite class of hairy cycles, the \emph{coronas} $C_{2t}\odot \lambda K_1$, see \cite{FH}.
%
%

\begin{rem}\label{rem}
Let $H=HC(n_1,m_1,\dots,n_t,m_t)$ with $t\geq 2$ be a bipartite hairy cycle. We will see $H$ as the hairy
cycle obtained from the caterpillar $\Gamma$ of Figure \ref{Sunxy} adding the edge $[x_1,y_t]$.
Clearly, the number of edges of $H$ is $e=2t+\sum^t_{i=1}n_i+\sum^t_{i=1}m_i$. Let $k=t+\sum^t_{i=1}n_i$, that is  $k=|B|$.
Let $e=d\cdot m$, $\Delta=\{1,2,\dots,e+d-1\}$, $\Delta'=\{m+1,2(m+1),\dots,(d-1)(m+1)\}$,
$S=(\Delta\setminus \Delta')\setminus\{c\}$ where $c\in \Delta\setminus\Delta'$.
Let $f$ be the standard $\alpha_S$-labeling of $\Gamma$. The missing integers
of $\Delta'$ and the removal of $c$ cause $d$ mv-labels in $f(A)\bigcup f(B)$,
let $d_A$ and $d_B$ be the number of mv-labels in $f(A)$ and
$f(B)$, respectively. Now $f(y_t)-f(x_1)=f(y_t)-(f(x_1)+1)+1=|B|+d_B=k+d_B$. Thus if $k+d_B=c$, in a natural way we can
extend $f$ to a  $d$-divisible $\alpha$-labeling of $H$, labeling the edge $[x_1,y_t]$ by $c$.
\end{rem}

\subsubsection{Odd (e-divisible) $\alpha$-labelings of hairy cycles}

\begin{thm}\label{thm:sunshine}
A hairy cycle admits an odd $\alpha$-labeling if and only if it is bipartite.
\end{thm}
\noindent
Proof. It is known that a hairy cycle with an odd $\alpha$-labeling
is necassarily bipartite.

Suppose now $H$ to be a bipartite hairy cycle.
We start from the above Remark \ref{rem} and
consider that $e=d\cdot m$, $d=e$ so $m=1$, hence
$\Delta\setminus\Delta'=\{1,3,5,\ldots,2e-1\}$.
Thus, if the removal of $c$ causes an
extra mv-label in $f(A)$ we will have $d_B=k-1$ and $f(y_t)-f(x_1)=2k-1$; if the removal of $c$ causes an
extra mv-label in $f(B)$ we will have $d_B=k+1$ and $f(y_t)-f(x_1)=2k+1$.
Therefore, it will be convenient to choose $c\in \{2k-1,2k+1\}$
(note that this choice is always possible since $\{2k-1,2k+1\}\subseteq \Delta\setminus\Delta'$).
In other words, the removal of $c$ does not affect the values of the labels
up to the $(k-1)$-th edge, the $k$-th edge will be labeled by $2k-1$ if we choose $c=2k+1$
and vice versa, and after the labeling continues in a natural way up to the end.
 There are four possible cases.\\
$\bullet$ Case (1)~If the $k$-th edge is a pendant edge from a vertex in $B$, the extra mv-label
caused by the removal of $c$ will be always
in $f(A)$, thus $f(y_t)-f(x_1)=2k-1$. Hence, choosing $c=2k-1$ and using
$2k+1$ as label of the $k$-th edge, we obtain an odd $\alpha$-labeling of $H$.\\
$\bullet$ Case (2)~If the $k$-th edge is a pendant edge from a vertex in $A$, the extra mv-label
caused by the removal of $c$ will be always
in $f(B)$, thus $f(y_t)-f(x_1)=2k+1$. Hence, choosing $c=2k+1$ and using
$2k-1$ as label of the $k$-th edge, we obtain an odd $\alpha$-labeling of $H$.\\
$\bullet$ Case (3)~If the $k$-th edge is an edge of the cycle of the form $[x_s,y_{s-1}]$, $c=2k-1$ implies
that the extra mv-label will be in $f(A)$, so $f(y_t)-f(x_1)=2k-1$, while $c=2k+1$ implies
that the extra mv-label will be in $f(B)$, so $f(y_t)-f(x_1)=2k+1$. In both cases, the choice of the value of $c$ results in appropriate.\\
$\bullet$ Case (4)~If the $k$-th edge is an edge of the cycle of the form $[x_s,y_{s}]$, $c=2k-1$ implies
that the extra mv-label will be in $f(B)$, so $f(y_t)-f(x_1)=2k+1$ while $c=2k+1$ implies
that the extra mv-label will be in $f(A)$, so $f(y_t)-f(x_1)=2k-1$. So, no choice results in
appropriate because the value we need to label $[x_1,y_t]$ has been already used.
To solve the problem, we have to distinguish two subcases: case $(4_1)$~$n_{s+1}\neq 0$, and case $(4_2)$~$n_t\neq 0 \bigwedge m_t=0$.
If we are able to define an odd $\alpha$-labeling in both previous subcases, we can do so in anycase.
In fact,
if the sequence $(n_1,m_1,\dots,n_t,m_t)$ does not contains zeros, we refer to the  Case $(4_1)$. If the sequence
$(n_1,m_1,\dots,n_t,m_t)$ contains at least one zero we refer to the  Case $(4_2)$ because we can always
choose another representation
 of $H$ so that it results in $n_t\neq 0 \ \bigwedge \ m_t=0$, $H$ being not a cycle.

Case $(4_1)$~~We choose $c=2k-1$ and applying $[O_2]$ we obtain another $\alpha_S$-labeling $g$
with $g(y_t)-g(x_1)=f(y_t)-f(x_1)-2=2k+1-2=2k-1$.
$[O_2]$ can be used because both (1) and (2) of its definition hold, in fact:\\
we know that $f(y_s)-f(x_s)=2k+1$ and $c=2k-1$ has been removed, so $f(x_s^{n_s})-f(x_s)=2k-3$. Hence
 $f(x_s^{n_s})+4=f(y_s)$ and (2) holds. Moreover, from $n_{s+1}\neq 0$ we have
$f(x_{s+1}^1)-f(x_{s+1})=f(y_s)-f(x_{s+1})+2$. Hence
$f(x_{s+1}^1)=f(y_s)+2$ and (1) holds.

Case $(4_2)$~~ We choose $c=2k-1$ and applying $[O_4]$ we obtain another $\alpha_S$-labeling $g$
with $g(y_t)-g(x_1)=f(y_t)-f(x_1)-2=2k-1$.
$[O_4]$ can be used because both (1) and (2) of its definition hold by the hypothesis.

\hfill $\Box$\\

\noindent
For an explicit definition of the odd $\alpha$-labeling described in the above theorem see the Appendix.

\begin{ex}
In Figure \ref{ExCase4b}(a) we have the graph $H=HC(3,3,0,0,3,6,0,1,$ $3,1)$. Clearly we have
$e=30$ and $k=14$. In particular, the $k$-th edge is $[x_3,y_3]$.
So, following the notation of the proof of Theorem \ref{thm:sunshine} Case (4), we have $s=3$.
Since $n_{s+1}=n_4=0$ and $m_5\neq 0$,
we have to rearrange the representation of $H$ as $HC(1,3,1,3,3,0,0,3,6,0)$.
It just so happens that the $k$-th edge is still
$[x_3,y_3]$, but now we are in the hypotheses of the case $(4_2)$, see Figure \ref{ExCase4b}(b).
The odd $\alpha$-labeling of $H$ constructed following the proof of the previous theorem is shown in Figure \ref{ExCase4b}(c).\\
Obviously, if after the rearrangement the $k$-th edge is not of the form $[x_s,y_s]$
we are in one of the other cases and we apply the corresponding construction.
\begin{figure}[H]
\begin{center}
\psfrag{x1}{$x_1$}
\psfrag{x2}{$x_2$}
\psfrag{x3}{$x_3$}
\psfrag{x4}{$x_4$}
\psfrag{x5}{$x_5$}
\psfrag{y1}{$y_1$}
\psfrag{y2}{$y_2$}
\psfrag{y3}{$y_3$}
\psfrag{y4}{$y_4$}
\psfrag{y5}{$y_5$}
\psfrag{0}{\tiny{$0$}}
\psfrag{1}{\tiny{$1$}}
\psfrag{2}{\tiny{$2$}}
\psfrag{3}{\tiny{$3$}}
\psfrag{4}{\tiny{$4$}}
\psfrag{5}{\tiny{$5$}}
\psfrag{6}{\tiny{$6$}}
\psfrag{7}{\tiny{$7$}}
\psfrag{8}{\tiny{$8$}}
\psfrag{9}{\tiny{$9$}}
\psfrag{10}{\tiny{$10$}}
\psfrag{11}{\tiny{$11$}}
\psfrag{12}{\tiny{$12$}}
\psfrag{13}{\tiny{$13$}}
\psfrag{14}{\tiny{$14$}}
\psfrag{15}{\tiny{$15$}}
\psfrag{16}{\tiny{$16$}}
\psfrag{17}{\tiny{$17$}}
\psfrag{18}{\tiny{$18$}}
\psfrag{19}{\tiny{$19$}}
\psfrag{20}{\tiny{$20$}}
\psfrag{21}{\tiny{$21$}}
\psfrag{22}{\tiny{$22$}}
\psfrag{23}{\tiny{$23$}}
\psfrag{24}{\tiny{$24$}}
\psfrag{25}{\tiny{$25$}}
\psfrag{26}{\tiny{$26$}}
\psfrag{27}{\tiny{$27$}}
\psfrag{28}{\tiny{$28$}}
\psfrag{29}{\tiny{$29$}}
\psfrag{30}{\tiny{$30$}}
\psfrag{31}{\tiny{$31$}}
\psfrag{32}{\tiny{$32$}}
\psfrag{33}{\tiny{$33$}}
\psfrag{34}{\tiny{$34$}}
\psfrag{35}{\tiny{$35$}}
\psfrag{36}{\tiny{$36$}}
\psfrag{37}{\tiny{$37$}}
\psfrag{38}{\tiny{$38$}}
\psfrag{39}{\tiny{$39$}}
\psfrag{40}{\tiny{$40$}}
\psfrag{41}{\tiny{$41$}}
\psfrag{43}{\tiny{$43$}}
\psfrag{45}{\tiny{$45$}}
\psfrag{47}{\tiny{$47$}}
\psfrag{49}{\tiny{$49$}}
\psfrag{51}{\tiny{$51$}}
\psfrag{53}{\tiny{$53$}}
\psfrag{55}{\tiny{$55$}}
\psfrag{57}{\tiny{$57$}}
\psfrag{59}{\tiny{$59$}}
\subfigure[]{\includegraphics[width=0.26\textwidth]{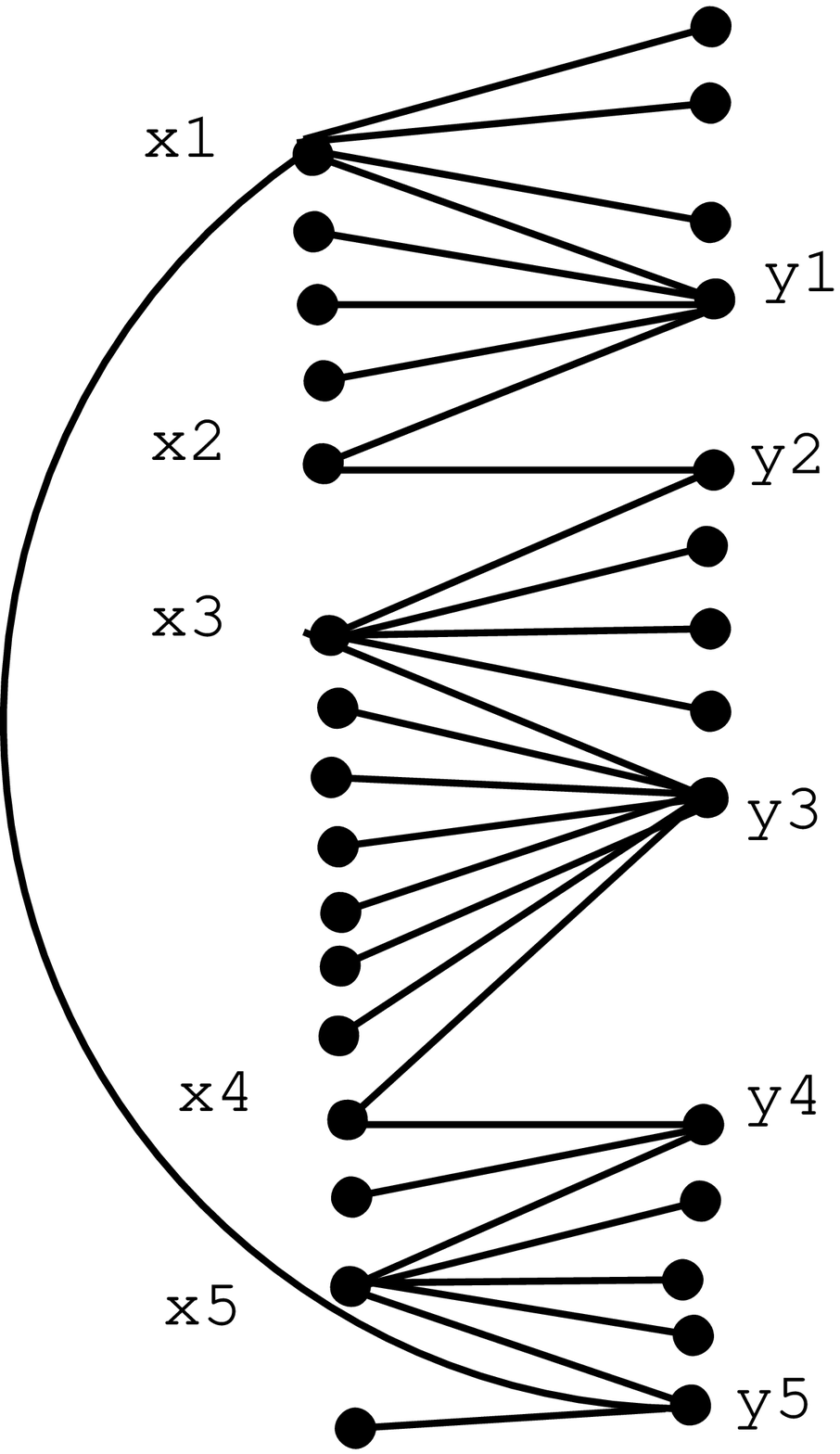}}\quad\quad\quad
\subfigure[]{\includegraphics[width=0.2\textwidth]{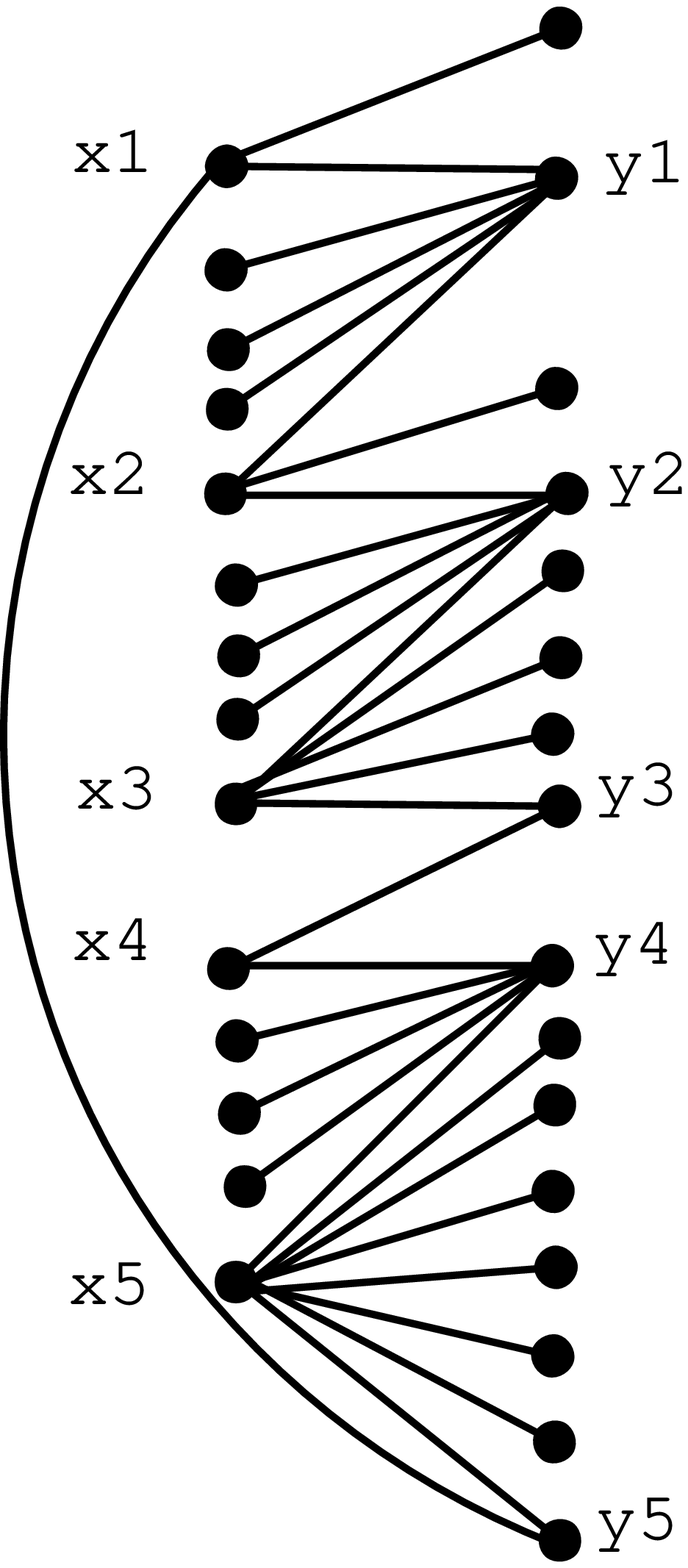}}\quad\quad\quad
\subfigure[]{\includegraphics[width=0.2\textwidth]{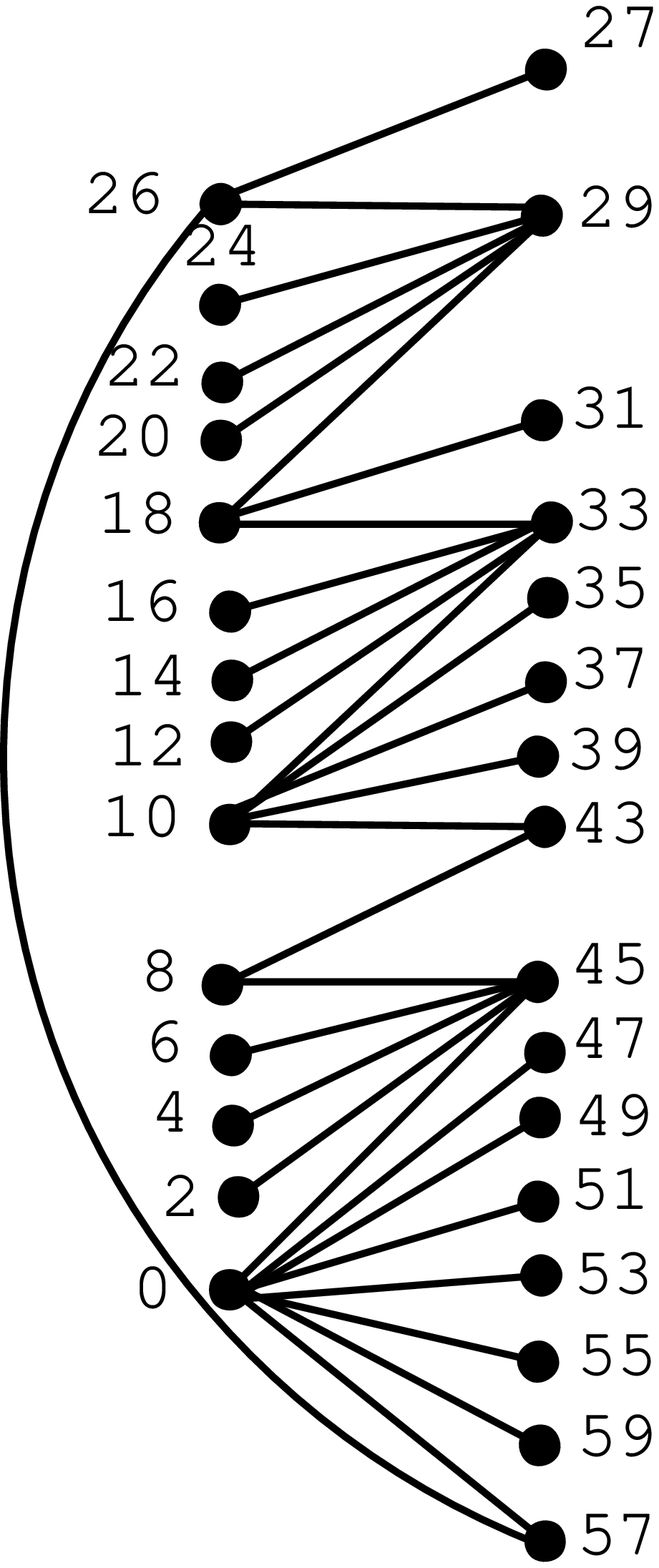}}
\caption{}\label{ExCase4b}
\end{center}
\end{figure}
\end{ex}
\noindent As a consequence of Theorems \ref{decomp} and \ref{thm:sunshine} we have:
\begin{thm}
Let $\G$ be a bipartite hairy cycle of size $e$. There exists a cyclic $\G$-decomposition of
$K_{2 \times 2en}$ for any positive integer $n$.
\end{thm}

\subsubsection{$d$-divisible $\alpha$-labelings of $C_{2t}\bigodot \lambda K_1$}
For convenience we will denote by $H(2t,\lambda)$
the hairy cycle $HC(\underbrace{\lambda,\lambda,\ldots,\lambda}_{2t})$
 with the cycle of length $2t$, $t\geq2$,
and with $\lambda$, $\lambda\geq1$, pendant edges through each vertex of the cycle.
From \cite{FH} we learn that such a graph is nothing but the
\emph{corona} of $C_{2t}$ with $\lambda K_1$, denoted by $C_{2t}\bigodot \lambda K_1$.
Obviously $H(2t,\lambda)$  has $2t(\lambda+1)$ edges.\\

\begin{thm}\label{H2tn}
The hairy cycle $H(2t,\lambda)$ admits a $d$-divisible $\alpha$-labeling
for any admissible value of $d$.
\end{thm}
\noindent
Proof. If $d=e$ the result follows from Theorem \ref{thm:sunshine}. So, from now on, we can assume $d\neq e$.
Let $\G$ be the caterpillar obtained deleting the edge $[x_1,y_t]$ from $H(2t,\lambda)$ and let $P$ be the path associated to $\G$.
We start again from the Remark \ref{rem} and we notice that now $e=d\cdot m=2t(\lambda+1)$, $k=e/2$ and the $k$-th edge is always an
edge of the path $P$: of the form $[x_r,y_r]$ if $t$ is odd, of the form $[x_r,y_{r-1}]$ if $t$ is even.

Also, if $d$ is odd the mv-labels due to the $d-1$ missing elements of $\Delta'$ are
equally distributed in $f(A)$ and $f(B)$,  because of the symmetry of the graph.
So, if the removal of $c$ causes an extra mv-label in $f(A)$, we will have
$d_B=\frac{d-1}{2}$  and $f(y_t)-f(x_1)=\frac{e+d-1}{2}$. If the removal of $c$ causes an extra mv-label in $f(B)$,
we will have $d_B=\frac{d+1}{2}$  and $f(y_t)-f(x_1)=\frac{e+d+1}{2}$. Thus, it will be convenient to choose
$c\in \{c_1,c_2 \}$ where $c_1=\frac{e+d-1}{2}$ and $c_2=\frac{e+d+1}{2}$.
Notice that $\{c_1,c_2\}\subseteq\Delta\setminus\Delta'$ as if by contradiction
$\frac{e+d\pm 1}{2}=\mu(m+1)$ then $(d-2\mu)(m+1)=\pm 1$, that is impossible.\\
If $d$ is even, because of the symmetry of the graph, the mv-labels due to the $d-2$ elements of
$\Delta'\setminus\{e/2\}$ are equally distributed in $f(A)$ and $f(B)$. So, if the removal of $c$ causes the two further
  mv-labels in $f(A)$, we will have $d_B=\frac{d-2}{2}$  and $f(y_t)-f(x_1)=\frac{e+d-2}{2}$ as well as
 if the removal of $c$ causes the two further mv-labels in $f(B)$, we will have $d_B=\frac{d+2}{2}$
 and $f(y_t)-f(x_1)=\frac{e+d+2}{2}$. Thus, it will be convenient to choose
$c\in \{c_1,c_2 \}$ where $c_1=\frac{e+d-2}{2}$ and $c_2=\frac{e+d+2}{2}$.
Notice that $\{c_1,c_2\}\subseteq\Delta\setminus\Delta'$ since $\frac{e+d\pm2}{2}=\frac{d}{2}(m+1)\pm 1 \notin \Delta'.$\\
 In other words, in both cases the removal of $c$ does not affect the values of the labels up to the $(k-1)$-th edge,
 the $k$-th edge will be labelled by $c_1$ if we have chosen $c=c_2$ and vice versa, and after the labeling continues,
 in ascending order, up to the end. We continue by distinguishing several cases.

\noindent $\bullet$~Case $(1)$~Let $t$ be even, that is the $k$-th edge is of the form $[x_r,y_{r-1}]$
with $r=\frac{t+2}{2}$.\\
Let $S=(\Delta\setminus \Delta')\setminus\{c_2\}$. Let $f$ be the standard $\alpha_S$-labeling of $\Gamma$.
The choice of $c=c_2$ implies that  the extra mv-labels  are always in $f(B)$, so $f(y_t)-f(x_1)=c_2$ and, in a natural way,
 $f$ can be extended to a $d$-divisible $\alpha$-labeling of $H(2t,\lambda)$, labeling the edge $[x_1,y_t]$ by $c_1$.\\

\noindent$\bullet$~Case $(2)$~Let $t$ be odd, that is the $k$-th edge is of the form $[x_r,y_{r}]$ with
$r=\frac{t+1}{2}$.\\
Let $S=(\Delta\setminus \Delta')\setminus\{c\}$ where $c$ can be chosen in $\{c_1,c_2\}$.
Let $f$ be the standard $\alpha_S$-labeling of $\Gamma$. The choice $c=c_1$ implies that the extra mv-labels  are always
in $f(B)$, so $f(y_t)-f(x_1)=c_2$. The choice $c=c_2$ implies that the extra mv-labels  are always in $f(A)$, so
$f(y_t)-f(x_1)=c_1$. In both cases the value we need  is not available. We can proceed by distinguishing two subcases. \\
\indent~Case $(2_1)$~Let $d$ be odd.
Now $\{c_1,c_2\}=\{\frac{e+d-1}{2},\frac{e+d+1}{2}\}$.
We choose $c=c_1$, consider the corresponding standard $\alpha_S$-labeling of $\Gamma$
 and apply $[O_1]$ to $f$ with $s=\frac{t+1}{2}$. This is possible because  (1) and (2) are satisfied, in fact
now we have $f(y_s)-f(x_s)=c_2$ and
$f(x_s^{n_s})-f(x_s)=c_1-1$ if $c_1-1\notin \Delta'$ or
$f(x_s^{n_s})-f(x_s)=c_1-2$ if $c_1-1\in \Delta'$, so
$f(y_s)-f(x_s^{n_s})=2$ or $3$ and (2) is true.
Also, if (1) was not true, there would be a mv-label between $f(y_s)$ and $f(x_{s+1}^{1})$ and, symmetrically,
there would be a mv-label between $f(x_{s})$ and $f(y_{s-1}^{m_{s-1}})$.
Thus $m$ should be a divisor of
$2\lambda+3$, an odd number, while $m$ is obviously even as now $d$ is odd.
Applying $[O_1]$ we obtain an $\alpha_S$-labeling $g$  of $\G$ with $g(y_t)-g(x_1)=f(y_t)-f(x_1)-1=c_2-1=c_1$.
So, we can
 extend $g$ to a $d$-divisible $\alpha$-labeling of $H(2t,\lambda)$,
 labeling the edge $[x_1,y_t]$ by $c_1$.\\
\indent~Case $(2_2)$~Let $d$ be even.
Now $\{c_1,c_2\}=\{\frac{e+d-2}{2},\frac{e+d+2}{2}\}$.
  We have to split the proof in several subcases.\\
\indent~Case $(2_{2_1})$~Let $\lambda\geq 2$ and $\lambda\not\equiv m-2,m-3$~(mod $m$).\\
Choose $c=c_1$ and let $f$ be the corresponding standard $\alpha_S$-labeling of $\Gamma$,
apply $[O_3]$ to $f$ with $s=\frac{t+1}{2}$. This is possible because  (1),(2) and (3) are satisfied, in fact
now we have $f(y_s)-f(x_s)=c_2$ and
$f(x_{s}^{n_s})-f(x_{s})=c_1-1$ if $c_1-1\notin \Delta'$ or
$f(x_s^{n_s})-f(x_s)=c_1-2$ if $c_1-1\in \Delta'$
so
$f(y_{s})-f(x_{s}^{n_{s}})=3$ or $4$ and (2) is true.
Also, if (1) was not true, there would be a mv-label either between $f(y_{s})$ and $f(x_{s+1}^{1})$
(excluded as $\lambda\not\equiv m-2$~(mod $m$)) or between $f(x_{s+1}^{1})$ and $f(x_{s+1}^{2})$
(excluded as $\lambda\not\equiv m-3$~(mod $m$)).
Applying $[O_3]$ we obtain an $\alpha_S$-labeling $g$  of $\G$ with
$g(y_t)-g(x_1)=f(y_t)-f(x_1)-2=c_2-2=c_1$.
So, we can
 extend $g$ to a $d$-divisible $\alpha$-labeling of $H(2t,\lambda)$,
 labeling the edge $[x_1,y_t]$ by $c_1$.\\
\indent~Case $(2_{2_2})$~Let $\lambda\geq2$ and $\lambda\equiv m-2$~(mod $m$).\\
Choose $c=c_1$ and let $f$ be the corresponding standard $\alpha_S$-labeling of $\Gamma$,
apply $[O_2]$ to $f$ with $s=\frac{t+1}{2}$. This is possible because  (1) and (2) are satisfied, in fact
now we have $f(y_{s})-f(x_{s})=c_2$ and
$f(x_{s}^{n_s})-f(x_{s})=c_1-1$ if $c_1-1\notin \Delta'$ or
$f(x_s^{n_s})-f(x_s)=c_1-2$ if $c_1-1\in \Delta'$
so
$f(y_{s})-f(x_{s}^{n_{s}})=3$ or $4$ and (2) is true.
Also, there is a mv-label between $f(y_{s})$ and $f(x_{s+1}^{1})$ as $\lambda\equiv m-2$~(mod $m$),
so $f(y_{s})-f(x_{s+1})+2=f(x_{s+1}^{1})-f(x_{s+1})$ and $(1)$ is true.
Applying $[O_2]$ we obtain an $\alpha_S$-labeling $g$  of $\G$ with
$g(y_t)-g(x_1)=f(y_t)-f(x_1)-2=c_2-2=c_1$.
So, we can
 extend $g$ to a $d$-divisible $\alpha$-labeling of $H(2t,\lambda)$,
 labeling the edge $[x_1,y_t]$ by $c_1$.\\
\indent~Case $(2_{2_3})$~Let $\lambda\equiv m-3$~(mod $m$) and \\
\indent$\bullet$~$\lambda\geq 4 \bigwedge m\geq 4\bigwedge t\geq 5 $ or
$\lambda\geq 7 \bigwedge m\geq 4\bigwedge t=3$.\\
Choose $c=c_1$, consider the corresponding standard $\alpha_S$-labeling of $\Gamma$, say $f$,
 and apply $[O_1]$ and $[O_5]_4$ to obtain an $\alpha_S$-labeling $g$ of $\G$.  We can apply  $[O_1]$  because
 $\lambda\equiv m-3$~(mod $m$),  and we can apply $[O_5]_4$ as $\lambda\geq 4 \bigwedge m\geq 4$.
 We can apply both $[O_1]$ and $[O_5]_4$, in any order, as
 $t\geq 5$ or $t=3 \bigwedge 7\leq \lambda$, and this ensures that there is no edge affected by both the
 operations. Thus we have
 $g(y_t)-g(x_1)=(f(y_t)-f(x_1)-1)-1=c_2-2=c_1$.
 So, we can
  extend $g$ to a $d$-divisible $\alpha$-labeling of $H(2t,\lambda)$,
  labeling the edge $[x_1,y_t]$ by $c_1$.\\
\indent$\bullet$~$\lambda\geq 4 \bigwedge m\geq 4\bigwedge t=3 \bigwedge \lambda<7$.\\
From our hypotheses we have

\begin{tabular}{l|l|l|l}
$4\leq \lambda < 7$&$\lambda+3$& $e=6(\lambda+1)$&$m\geq 4$ divides $\lambda+3$ and $e$\\
\hline
4&7&30& there is no value \\
5&8&36& 4 unacceptable, as $d$ is even \\
6&9&42& there is no value \\
\hline
\end{tabular}

\smallskip

$\bullet$~$m=2$, $t\geq 3$.\\
Choose $c=c_1$, consider the corresponding standard $\alpha_S$-labeling of $\Gamma$, say $f$,
 and apply $[O_1]$ with $s=\frac{t+1}{2}$ and $[O_5]$ with $s=1$ to obtain an $\alpha_S$-labeling $g$ of $\G$.
 We can apply $[O_1]$  because  $\lambda\equiv m-3$~(mod $m$),  and we can apply  $[O_5]$ as $m$ divides $2(\lambda+1)$.
 We can apply both $[O_1]$ and $[O_5]$, in any order, as $\lambda\geq 2$ and this ensures that there is no edge
 affected by both the operations.\\
Thus we have  $g(y_t)-g(x_1)=(f(y_t)-f(x_1)-1)-1=c_2-2=c_1$.
So we can extend $g$ to a $d$-divisible $\alpha$-labeling of $H(2t,\lambda)$,
labeling the edge $[x_1,y_t]$ by $c_1$.\\
\indent$\bullet$~$m=3$, $t\geq 3$.\\
Choose $c=c_1$ and consider the corresponding standard $\alpha_S$-labeling of $\Gamma$, say $f$.
We can apply $[O_1]$ with $s=\frac{t+1}{2}$ because  $\lambda\equiv m-3$~(mod $m$),  and we can apply  $[O_1]$ with $s=1$ as
$m$ divides $\lambda$.  We can apply twice $[O_1]$ with $s=1$ and $s=\frac{t+1}{2}$, in any order, as $\lambda\geq 2$ and this ensures that there is no edge
affected by both the operations.\\
Thus we have  $g(y_t)-g(x_1)=(f(y_t)-f(x_1)-1)-1=c_2-2=c_1$.
So, we can extend $g$ to a $d$-divisible $\alpha$-labeling of $H(2t,\lambda)$,
labeling the edge $[x_1,y_t]$ by $c_1$.\\
 \indent$\bullet$~$m\geq 4 \bigwedge \lambda=2.$\\
 It happens only when $m=5$ e $\lambda=2$. Then $e=6t$ implies that $t$ have to be an odd multiple of 5.
 We choose $c=c_1$ and consider the corresponding standard $\alpha_S$-labeling of $\Gamma$, say $f$.
 It easy to see that we can apply both $[O_1]$ with $s=\frac{t+1}{2}$ and $[O_5]$ with $s=\frac{t+3}{2}$, in any order, being sure that there is no edge
 affected by both the operations.\\
 Thus we have  $g(y_t)-g(x_1)=(f(y_t)-f(x_1)-1)-1=c_2-2=c_1$.
 So, we can extend
 $g$ to a $d$-divisible $\alpha$-labeling of $H(2t,\lambda)$,
 labeling the edge $[x_1,y_t]$ by $c_1$.\\
 \indent$\bullet$~$m\geq 4 \bigwedge \lambda=3.$\\
 It happens only when $m=6$ e $\lambda=3$. Then $e=8t$ implies that $t$ have to be an odd multiple of 3.
If $t=3$,  we are able to give directly a $d$-divisible
 $\alpha$-labeling of $H(2t,\lambda)$.
 If $t\geq 5$, we choose $c=c_1$ and consider the corresponding standard $\alpha_S$-labeling of $\Gamma$, say $f$.
 It is easy to see that we can apply both $[O_1]$ with $s=\frac{t+1}{2}$ and $[O_5]$ with $s=3$,
 in any order, being sure that there is no edge
 affected by both the operations.\\
 Thus we have  $g(y_t)-g(x_1)=(f(y_t)-f(x_1)-1)-1=c_2-2=c_1$.
 So, we can
 extend $g$ to a $d$-divisible $\alpha$-labeling of $H(2t,\lambda)$,
 labeling the edge $[x_1,y_t]$ by $c_1$.\\
 If $t=3$,  we are able to give directly a $d$-divisible
 $\alpha$-labeling of $H(2t,\lambda)$.

 Case $(2_{2_4})$ Let $\lambda=1$. Also in this case the foregoing construction can be applied,
 by distinguishing the cases $m=0,1,2,3 ($mod $4)$ and applying in a suitable way the $[O_i]$'s.
 Here, in Appendix, is given an explicit $2$-divisible $\alpha$-labeling
 of $H(2t,1)$ for  any $t$ odd.\hfill $\Box$

\begin{ex}
Here we show a $6$-divisible $\alpha$-labeling of the hairy cycle $HC(10,2)$. Since $t=5$, $d=6$,
$m=5$ and $\lambda=2$, we are in the Case $(2_{2_3})$ of Theorem \ref{H2tn}.   In Figure \ref{H10,2} (a)
we have the standard $\alpha_S$-labeling $f$ of the caterpillar $\G$ obtained by $HC(10,2)$ deleting
the edge $[x_1,y_{5}]$ with $S=(\Delta\setminus\Delta')\setminus\{17\}$ where $\Delta=\{1,2,\ldots,35\}$ and $\Delta'=\{6,12,18,24,30\}$.
Now $f(y_5)-f(x_1)=19$ which is not available, so we cannot extend $f$ to a $6$-divisible $\alpha$-labeling of $HC(10,2)$.
It is easy to see that it is possible to apply $[O_1]$ with $s=3$, the $\alpha_S$-labeling $g$ so obtained is shown
in Figure \ref{H10,2} (b).
Now $g(y_5)-g(x_1)=18$ which is again not available as $18 \in \Delta'$.
Then we can apply $[O_5]$ with $s=4$. After this we obtain
the $\alpha_S$-labeling $\tilde{g}$ of $\G$ shown in Figure \ref{H10,2} (c).
Finally $\tilde{g}(y_5)-\tilde{g}(x_1)=17$. Now we can extend $\tilde{g}$ to a $6$-divisible $\alpha$-labeling
of $HC(10,2)$, labeling the edge $[x_1,y_5]$ by $17$.
\begin{figure}[H]
\begin{center}
\psfrag{0}{\small{$0$}}
\psfrag{1}{\small{$1$}}
\psfrag{2}{\small{$2$}}
\psfrag{3}{\small{$3$}}
\psfrag{4}{\small{$4$}}
\psfrag{5}{\small{$5$}}
\psfrag{6}{\small{$6$}}
\psfrag{7}{\small{$7$}}
\psfrag{8}{\small{$8$}}
\psfrag{9}{\small{$9$}}
\psfrag{10}{\small{$10$}}
\psfrag{11}{\small{$11$}}
\psfrag{12}{\small{$12$}}
\psfrag{13}{\small{$13$}}
\psfrag{14}{\small{$14$}}
\psfrag{15}{\small{$15$}}
\psfrag{16}{\small{$16$}}
\psfrag{17}{\small{$17$}}
\psfrag{18}{\small{$18$}}
\psfrag{19}{\small{$19$}}
\psfrag{20}{\small{$20$}}
\psfrag{21}{\small{$21$}}
\psfrag{22}{\small{$22$}}
\psfrag{23}{\small{$23$}}
\psfrag{24}{\small{$24$}}
\psfrag{25}{\small{$25$}}
\psfrag{26}{\small{$26$}}
\psfrag{27}{\small{$27$}}
\psfrag{28}{\small{$28$}}
\psfrag{29}{\small{$29$}}
\psfrag{30}{\small{$30$}}
\psfrag{31}{\small{$31$}}
\psfrag{32}{\small{$32$}}
\psfrag{33}{\small{$33$}}
\psfrag{34}{\small{$34$}}
\psfrag{35}{\small{$35$}}
\psfrag{s1}{\tiny{$1$}}
\psfrag{s2}{\tiny{$2$}}
\psfrag{s3}{\tiny{$3$}}
\psfrag{s4}{\tiny{$4$}}
\psfrag{s5}{\tiny{$5$}}
\psfrag{s6}{\tiny{$6$}}
\psfrag{s7}{\tiny{$7$}}
\psfrag{s8}{\tiny{$8$}}
\psfrag{s9}{\tiny{$9$}}
\psfrag{s10}{\tiny{$10$}}
\psfrag{s11}{\tiny{$11$}}
\psfrag{s12}{\tiny{$12$}}
\psfrag{s13}{\tiny{$13$}}
\psfrag{s14}{\tiny{$14$}}
\psfrag{s15}{\tiny{$15$}}
\psfrag{s16}{\tiny{$16$}}
\psfrag{s17}{\tiny{$17$}}
\psfrag{s18}{\tiny{$18$}}
\psfrag{s19}{\tiny{$19$}}
\psfrag{s20}{\tiny{$20$}}
\psfrag{s21}{\tiny{$21$}}
\psfrag{s22}{\tiny{$22$}}
\psfrag{s23}{\tiny{$23$}}
\psfrag{s24}{\tiny{$24$}}
\psfrag{s25}{\tiny{$25$}}
\psfrag{s26}{\tiny{$26$}}
\psfrag{s27}{\tiny{$27$}}
\psfrag{s28}{\tiny{$28$}}
\psfrag{s29}{\tiny{$29$}}
\psfrag{s30}{\tiny{$30$}}
\psfrag{s31}{\tiny{$31$}}
\psfrag{s32}{\tiny{$32$}}
\psfrag{s33}{\tiny{$33$}}
\psfrag{s34}{\tiny{$34$}}
\psfrag{s35}{\tiny{$35$}}
\psfrag{b6}{\small{$\mathbf{6}$}}
\psfrag{b7}{\small{$\mathbf{7}$}}
\psfrag{b8}{\small{$\mathbf{8}$}}
\psfrag{b9}{\small{$\mathbf{9}$}}
\psfrag{b10}{\small{$\mathbf{10}$}}
\psfrag{b11}{\small{$\mathbf{11}$}}
\psfrag{b12}{\small{$\mathbf{12}$}}
\psfrag{b13}{\small{$\mathbf{13}$}}
\psfrag{b14}{\small{$\mathbf{14}$}}
\psfrag{b15}{\small{$\mathbf{15}$}}
\psfrag{b16}{\small{$\mathbf{16}$}}
\psfrag{b17}{\small{$\mathbf{17}$}}
\psfrag{b18}{\small{$\mathbf{18}$}}
\psfrag{b19}{\small{$\mathbf{19}$}}
\psfrag{b20}{\small{$\mathbf{20}$}}
\psfrag{b21}{\small{$\mathbf{21}$}}
\psfrag{b22}{\small{$\mathbf{22}$}}
\psfrag{b23}{\small{$\mathbf{23}$}}
\psfrag{b24}{\small{$\mathbf{24}$}}
\psfrag{b25}{\small{$\mathbf{25}$}}
\psfrag{b26}{\small{$\mathbf{26}$}}
\psfrag{b27}{\small{$\mathbf{27}$}}
\psfrag{b28}{\small{$\mathbf{28}$}}
\psfrag{b29}{\small{$\mathbf{29}$}}
\psfrag{b32}{\small{$\mathbf{32}$}}
\psfrag{bs22}{\tiny{$\mathbf{22}$}}
\psfrag{bs23}{\tiny{$\mathbf{23}$}}
\psfrag{bs25}{\tiny{$\mathbf{25}$}}
\psfrag{bs26}{\tiny{$\mathbf{26}$}}
\subfigure[]{
$\begin{array}{c}
\includegraphics[width=0.16\textwidth]{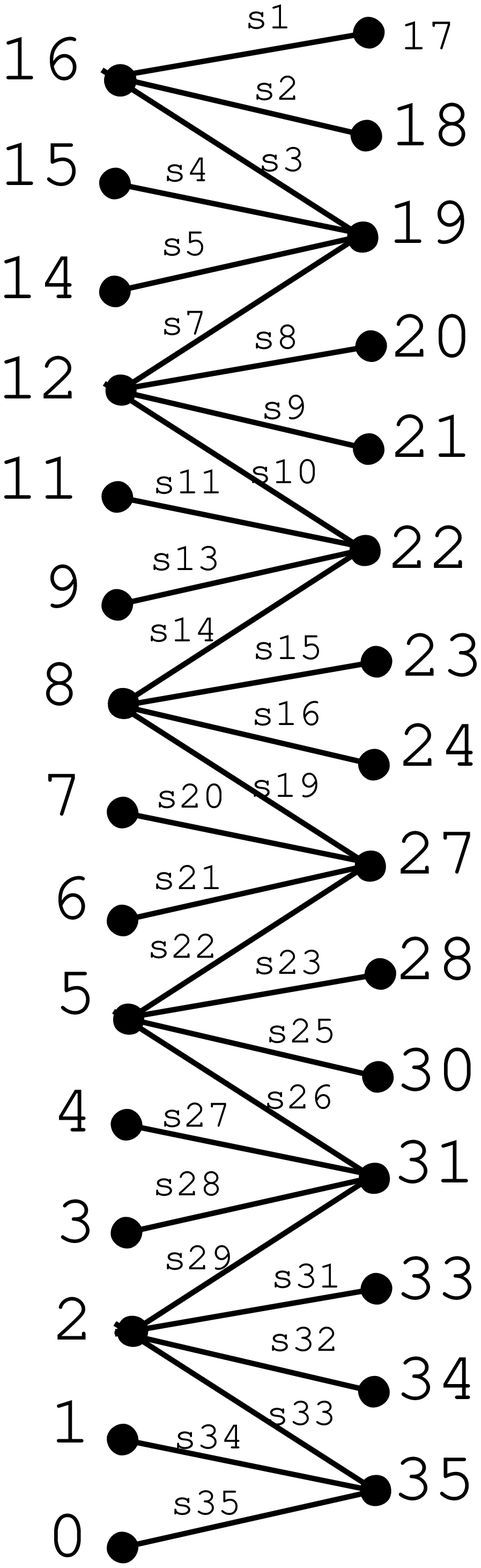}
\end{array}$}\quad\quad\quad
\subfigure[]{
$\begin{array}{c}
\includegraphics[width=0.19\textwidth]{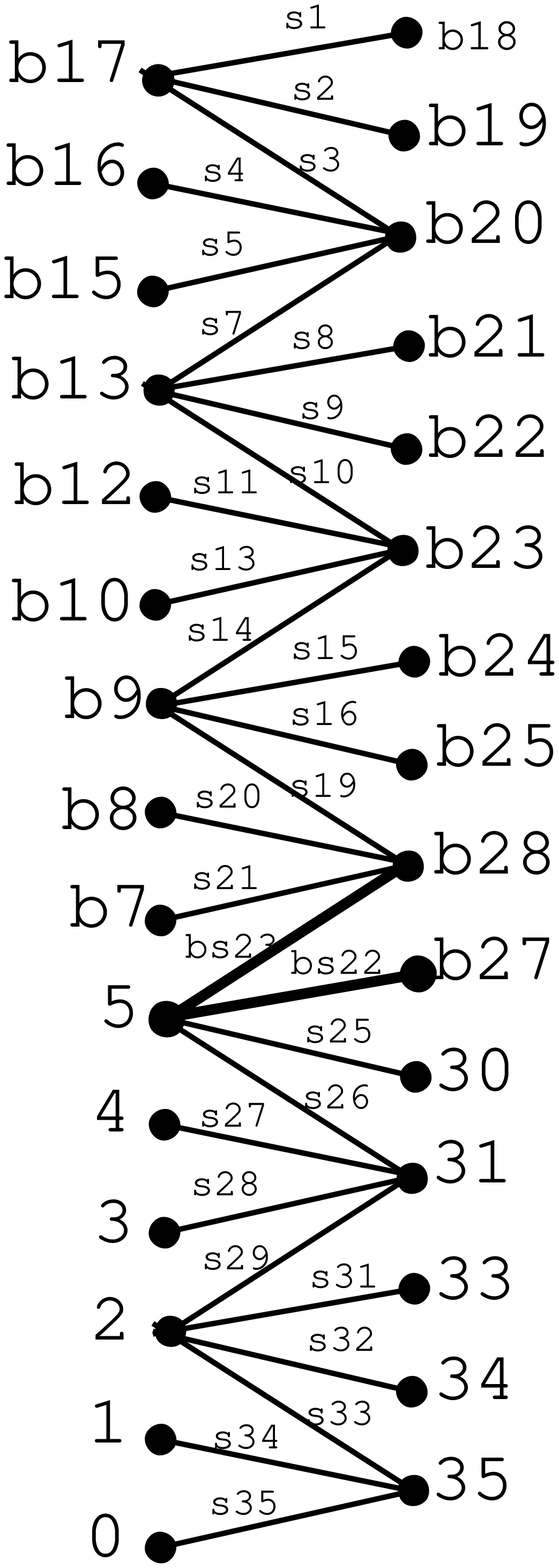}
\end{array}$}\quad\quad\quad
\subfigure[]{
$\begin{array}{c}
\includegraphics[width=0.19\textwidth]{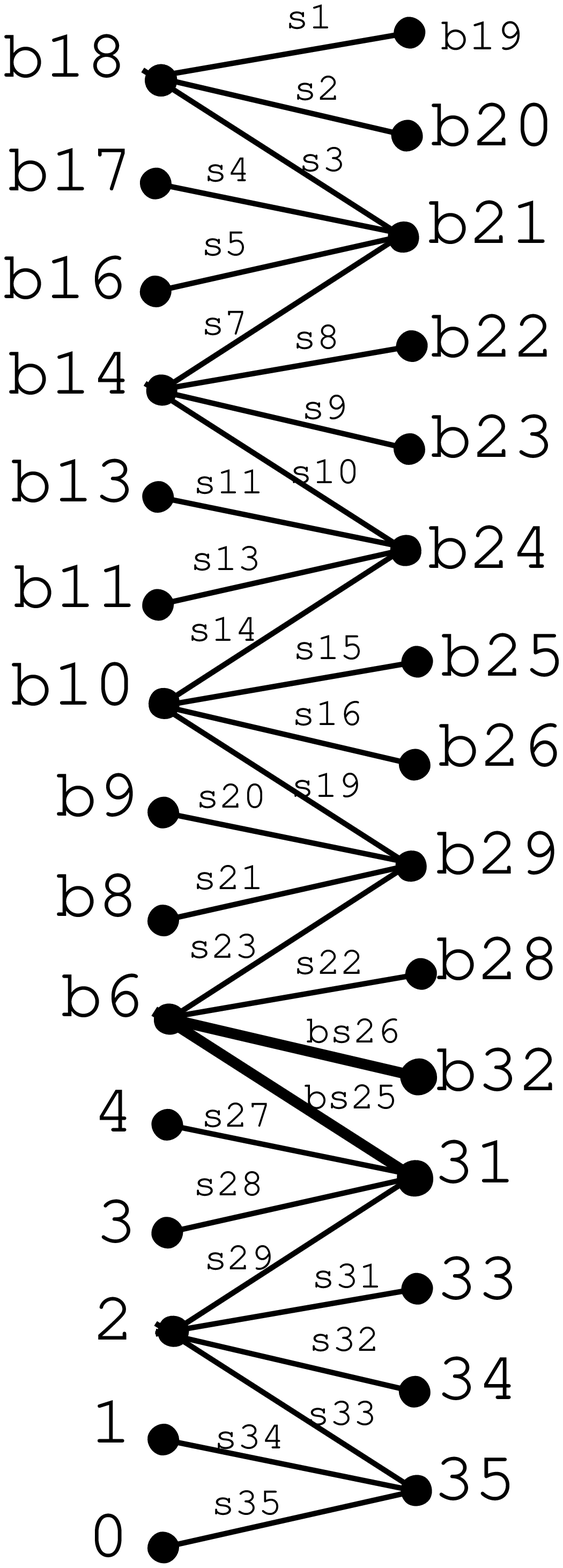}
\end{array}$}
\caption{}\label{H10,2}
\end{center}
\end{figure}
\end{ex}
\noindent
The following is an immediate consequence of Theorems \ref{decomp} and \ref{H2tn}.
\begin{thm}
There exists a cyclic $H(2t,\lambda)$-decomposition of
$K_{\left (\frac{2t(\lambda+1)}{d}+1\right )\times 2dn}$ for any admissible  $d$ and for any positive integer $n$.
\end{thm}

\subsection{Cycles}
As usual, we will denote  the cycle on $k$ vertices by $C_k$, $k\geq3$. It is obvious that
$C_k$ is a graph of size $k$ and that it is bipartite if and only if
$k$ is even.\\
In \cite{R} Rosa proved that $C_k$ has an $\alpha$-labeling if and only if
$k\equiv$0(mod 4). In \cite{APArs}, the second author proved that
$C_{4k}$ admits a $2$-divisible and a $4$-divisible $\alpha$-labeling for any positive integer $k$.
Here, generalizing this last result, we prove that $C_{4k}$
admits a $d$-divisible $\alpha$-labeling for any divisor $d$ of $4k$.

\begin{thm}\label{cycle}
For any positive integer $k$, the cycle $C_{4k}$ admits a $d$-divisible $\alpha$-labeling
for any admissible value of $d$.
\end{thm}
\noindent
Proof. Consider the cycle $C_{4k}$ as a bipartite graph as follows:
\begin{figure}[H]
\begin{center}
\psfrag{x1}{$x_1$}
\psfrag{x2}{$x_2$}
\psfrag{x3}{$x_3$}
\psfrag{x2k}{$x_{2k}$}
\psfrag{y1}{$y_1$}
\psfrag{y2}{$y_2$}
\psfrag{y3}{$y_3$}
\psfrag{y2k}{$y_{2k}$}
\psfrag{y2k-1}{$y_{2k-1}$}
\includegraphics[width=0.6\textwidth]{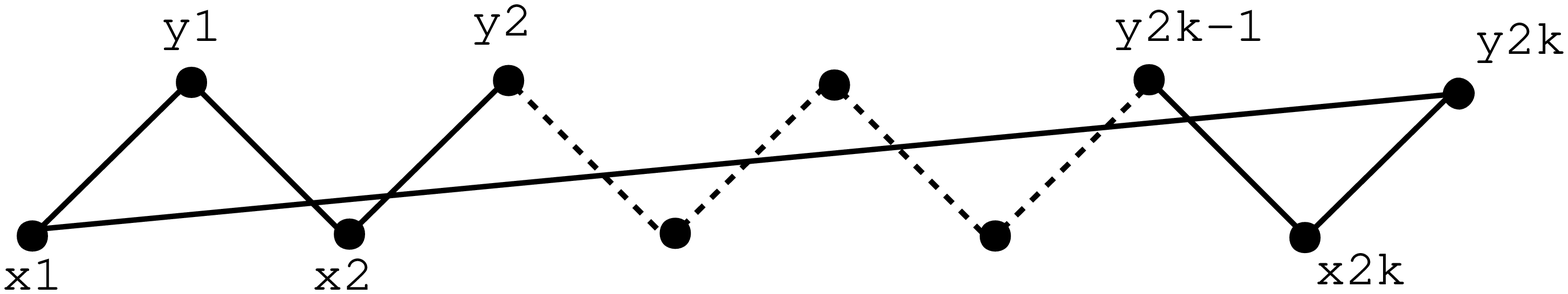}
\end{center}
\end{figure}
\noindent
and set $A=\{x_1,x_2,\dots,x_{2k}\}$ and $B=\{y_1,y_2,\dots,y_{2k}\}$ the two bipartite sets.
Let $\G$ be the caterpillar obtained form $C_{4k}$ deleting the edge $[x_1,y_{2k}]$.
Let $4k=d\cdot m$,
$\Delta=\{1,2,\dots,4k+d-1\}$  and $\Delta'=\{m+1,2(m+1),\dots,(d-1)(m+1)\}$.
Choose an element in $\Delta\setminus \Delta'$,
say $c$. Let $f$ be the standard $\alpha_S$-labeling of $\G$
where $S=(\Delta\setminus\Delta')\setminus\{c\}$.
In order to show that $f$ can be naturally extended to a $d$-divisible $\alpha$-labeling of $C_{4k}$ it remains
to prove that we can choose the element $c$ so that $c=f(y_{2k})-f(x_{1})$, see Remark \ref{rem:c}.\\
Let $d_x$ and $d_y$ denote the numbers of the mv-labels in $f(A)$ and $f(B)$, respectively.
Obviously, the condition $c=f(y_{2k})-f(x_1)$ becomes $c=|B|+d_y=2k+d_y$.
In what follows, we are able to determine $d_y$ and, consequently, $c$. \\
If $m=4k/d$ is even, the elements of $\Delta'$ are odd and even, alternatively, so $c$ lies between two elements of
different parity.
 The deletion of any element of $\Delta'$ smaller than $c$ produces a mv-label in $f(B)$, while deleting an element of
 $\Delta'$ greater than $c$ gives a mv-label in $f(A)$. Thus we have necessarily $d_y(m+1)<c<(d_y+1)(m+1)$, where $c=d_y+2k$ and
$2k=md/2$. With a simple calculation we obtain $d_ym<md/2\leq(d_y+1)m$, from which $2d_y<d\leq 2d_y+2$.
So,  $d$ even implies $d_y=(d-2)/2$ as well as $d$ odd implies $d_y=(d-1)/2$.\\
If $m=4k/d$ is odd, all the elements of $\Delta'$ are even.  If $c$ was even too, we would have $d_x=d_y=d/2$,
 so $c=2k+d/2=d/2(m+1)\in \Delta'$, but $c$ must not belong to $\Delta'$. Thus $c$ must be odd and we have $d_x-d_y=\pm 2$.
 In addition, we know that $d_x+d_y=d$, so $d_y=(d\mp 2)/2$.
 \hfill $\Box$\\
\begin{rem}
It is known that if $f$ is a ($d$-divisible) $\alpha$-labeling of a bipartite graph $\G$ of size $e$,
the function $g:V(\G)\rightarrow\{0,1,\ldots,e\}$,
defined by $g(x)=e-f(x)$, $\forall x \in V(\G)$ is again a ($d$-divisible) $\alpha$-labeling of $\G$.
We point out that if $d=1$, called $f$ the $\alpha$-labeling constructed in the above theorem and $g$
the classical $\alpha$-labeling given by A. Rosa in
\cite{R} it results $g(x)=e-f(x)$, $\forall x \in V(\G)$.
Also, when $d=2,4$ the same relation holds between the $d$-divisible $\alpha$-labeling constructed in the above theorem
 and that given by the second author in \cite{APArs}.
\end{rem}
\begin{ex}
In Figure \ref{c24} we show the $d$-divisible $\alpha$-labelings of $C_{24}$ described in Theorem \ref{cycle}
 for $d=3$ and $d=8$.\\
If $d=3$, we have $m+1=9$ and $c=13$. If $d=8$ it results $m+1=4$ and  we can choose $c=15$ or $c=17$ (as in Figure \ref{c24}).
\begin{figure}[H]
\begin{center}
\psfrag{1}{$1$}
\psfrag{2}{$2$}
\psfrag{3}{$3$}
\psfrag{0}{$0$}
\psfrag{4}{$4$}
\psfrag{5}{$5$}
\psfrag{6}{$6$}
\psfrag{7}{$7$}
\psfrag{8}{$8$}
\psfrag{9}{$9$}
\psfrag{10}{$10$}
\psfrag{11}{$11$}
\psfrag{13}{$13$}
\psfrag{12}{$12$}
\psfrag{14}{$14$}
\psfrag{15}{$15$}
\psfrag{16}{$16$}
\psfrag{17}{$17$}
\psfrag{18}{$18$}
\psfrag{19}{$19$}
\psfrag{20}{$20$}
\psfrag{21}{$21$}
\psfrag{22}{$22$}
\psfrag{23}{$23$}
\psfrag{24}{$24$}
\psfrag{25}{$25$}
\psfrag{26}{$26$}
\psfrag{27}{$27$}
\psfrag{28}{$28$}
\psfrag{29}{$29$}
\psfrag{30}{$30$}
\psfrag{31}{$31$}
\psfrag{1s}{{\tiny $1$}}
\psfrag{2s}{{\tiny $2$}}
\psfrag{6s}{{\tiny $6$}}
\psfrag{3s}{{\tiny $3$}}
\psfrag{4s}{{\tiny $4$}}
\psfrag{5s}{{\tiny $5$}}
\psfrag{7s}{{\tiny $7$}}
\psfrag{8s}{{\tiny $8$}}
\psfrag{9s}{{\tiny $9$}}
\psfrag{10s}{{\tiny $10$}}
\psfrag{11s}{{\tiny $11$}}
\psfrag{12s}{{\tiny $12$}}
\psfrag{14s}{{\tiny $14$}}
\psfrag{13s}{{\tiny $13$}}
\psfrag{15s}{{\tiny $15$}}
\psfrag{16s}{{\tiny $16$}}
\psfrag{17s}{{\tiny $17$}}
\psfrag{18s}{{\tiny $18$}}
\psfrag{19s}{{\tiny $19$}}
\psfrag{20s}{{\tiny $20$}}
\psfrag{21s}{{\tiny $21$}}
\psfrag{22s}{{\tiny $22$}}
\psfrag{23s}{{\tiny $23$}}
\psfrag{24s}{{\tiny $24$}}
\psfrag{25s}{{\tiny $25$}}
\psfrag{26s}{{\tiny $26$}}
\psfrag{27s}{{\tiny $27$}}
\psfrag{28s}{{\tiny $28$}}
\psfrag{29s}{{\tiny $29$}}
\psfrag{30s}{{\tiny $30$}}
\psfrag{31s}{{\tiny $31$}}
\includegraphics[width=0.9\textwidth]{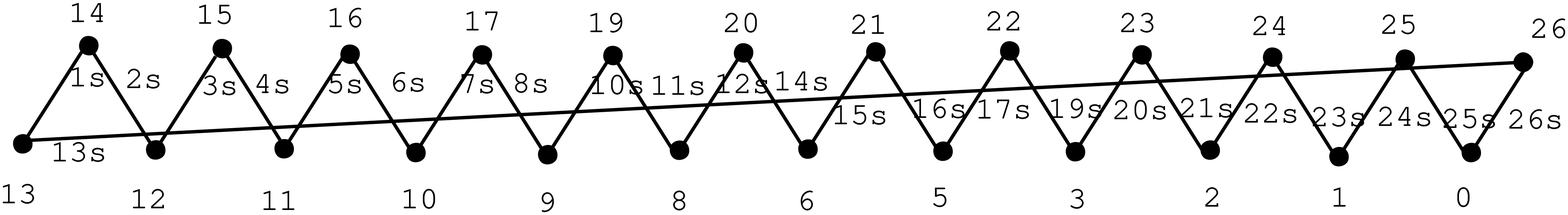}\\
\vspace{0.6cm}
\includegraphics[width=0.9\textwidth]{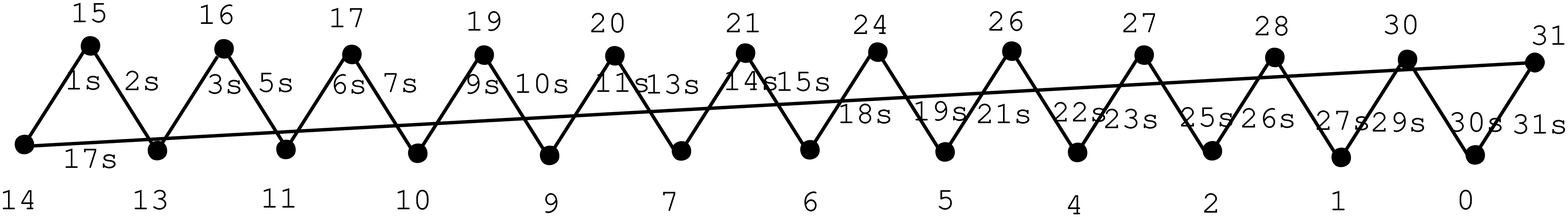}
\caption{}
\label{c24}
\end{center}
\end{figure}
\end{ex}
\noindent The following result immediately follows from Theorems \ref{decomp} and \ref{cycle}.
\begin{thm}
There exists a cyclic $C_{4k}$-decomposition of $K_{(\frac{4k}{d}+1)\times 2dn}$ for any
positive integers $k,n$ and any divisor $d$ of $4k$.
\end{thm}

\section*{Appendix}
Here we give an explicit representation of the odd $\alpha$-labeling
$f$ of a bipartite hairy cycle $\G$ of size $e$ constructed in Theorem
\ref{thm:sunshine}.
We give the definition of the function
 $f:V(\G)\rightarrow \{0,1,\ldots,2e-1\}$ in each of the possible cases.\\

\noindent
Case (1): the $k$-th edge is a pendant edge from a vertex in $B$, say $[y_s,y_s^j]$.
$$f(x_r)=\begin{cases}
2\sum_{\ell=r}^{t}m_\ell+2(t-r)+2\quad & \ r=1,\ldots,s\\
2\sum_{\ell=r}^{t}m_\ell+2(t-r)\quad & \ r=s+1,\ldots,t
\end{cases}$$
$$f(y_r)=f(x_1)+2\sum_{\ell=1}^{r}n_\ell+2r-1\quad \ r=1,\ldots,t$$
$$f(x_r^i)=\begin{cases}
f(x_1)+2i-1\quad & \ r=1,\ i=1,\ldots,n_1\\
f(y_{r-1})+2i\quad & \ r=2,\ldots,t, \ i=1,\ldots,n_r
\end{cases}$$
$$f(y_r^i)=\begin{cases}f(x_r)-2i\quad & \ r\neq s,\ i=1,\ldots,m_r\\
& \ r=s,\ i=1,\ldots,j-1\\
f(x_r)-2(i+1)\quad & \ r=s,\ i=j,\ldots,m_s.
\end{cases}$$

\noindent
Case (2): the $k$-th edge is a pendant edge from a vertex in $A$, say $[x_s,x_s^j]$.
$$f(x_r)=2\sum_{\ell=r}^{t}m_\ell+2(t-r)\quad \ r=1,\ldots,t$$
$$f(y_r^i)=f(x_r)-2i\quad \ r=1,\ldots,t,\ i=1,\ldots,m_r$$
$$f(y_r)=\begin{cases}
f(x_1)+2\sum_{\ell=1}^{r}n_\ell+2r-1\quad & \ r=1,\ldots,s-1\\
f(x_1)+2\sum_{\ell=1}^{r}n_\ell+2r+1\quad & \ r=s,\ldots,t
\end{cases}$$
$$f(x_r^i)=\begin{cases}
f(x_1)+2i-1\quad & \ r=1,\ i=1,\ldots,n_1\\
f(y_{r-1})+2i\quad & \ r\neq1,s,\ i=1,\ldots,n_r\\
\ & \ r=s,\ i=1,\ldots,j\\
f(y_{r-1})+2i+2\quad &\ r=s,\ i=j+1,\ldots,n_s,\ \textrm{if}\ j\neq n_s.
\end{cases}$$

\noindent
Case (3): the $k$-th edge is an edge of the cycle of the form $[x_s,y_{s-1}]$.
$$f(x_r)=\begin{cases}
2\sum_{\ell=r}^{t}m_\ell+2(t-r)+2\quad & \ r=1,\ldots,s-1\\
2\sum_{\ell=r}^{t}m_\ell+2(t-r)\quad & \ r=s,\ldots,t
\end{cases}$$
$$f(y_r)=f(x_1)+2\sum_{\ell=1}^{r}n_\ell+2r-1\quad \ r=1,\ldots,t$$
$$f(x_r^i)=\begin{cases}
f(x_1)+2i-1\quad & \ r=1,\ i=1,\ldots,n_1\\
f(y_{r-1})+2i\quad & \ r=2,\ldots,t, \ i=1,\ldots,n_r
\end{cases}$$
$$f(y_r^i)=f(x_r)-2i\quad \ r=1,\ldots,t,\ i=1,\ldots,m_r.$$

\noindent
Case $(4)$: the $k$-th edge is an edge of the cycle of the form $[x_s,y_{s}]$.\\
Case $(4_1):n_{s+1}\neq0$.
$$f(x_r)=\begin{cases}
2\sum_{\ell=r}^{t}m_\ell+2(t-r)+2\quad & \ r=1,\ldots,s\\
2\sum_{\ell=r}^{t}m_\ell+2(t-r)\quad & \ r=s+1,\ldots,t
\end{cases}$$
$$f(y_r^i)=f(x_r)-2i\quad  \ r=1,\ldots,t,\ i=1,\ldots,m_r$$
$$f(y_r)=\begin{cases}
f(x_1)+2\sum_{\ell=1}^{r}n_\ell+2r-1\quad & \ r\neq s\\
f(x_1)+2\sum_{\ell=1}^{s}n_\ell+2s+1\quad & \ r=s
\end{cases}$$
$$f(x_r^i)=\begin{cases}
f(x_1)+2i-1\quad & \ r=1, \ i=1,\ldots,n_1\\
f(y_{r-1})+2i\quad & \ r=2,\ldots,s,s+2,\ldots,t, \ i=1,\ldots,n_r\\
f(y_s)-2  & \ r=s+1, \ i=1\\
f(y_s)+2(i-1) & \  r=s+1, \ i=2,\ldots,n_{s+1}.
\end{cases}$$

\noindent
Case $(4_2):n_t\neq0$, $m_t=0$.
$$f(x_r)=2\sum_{\ell=r}^{t}m_\ell+2(t-r)\quad  \ \ \ \ \ r=1,\ldots,t$$
$$f(y_r^i)=f(x_r)-2i\quad  \ r=1,\ldots,t-1, \ i=1,\ldots,m_r$$
$$f(x_r^i)=\begin{cases}
f(x_1)+2i-1\quad & \ r=1, \ i=1,\ldots,n_1\\
f(y_{r-1})+2i-1\quad & \ r=2,\ldots,t-1,\ i=1,\ldots,n_r\\
& \ r=t,\ i=1,\ldots,{n_t-1}\\
2e-1\quad  & \ r=t,\ i=n_t
\end{cases}$$
$$f(y_r)=\begin{cases}
f(x_1)+2\sum_{\ell=1}^{r}n_\ell+2r-1\quad & \ r=1,\ldots,s-1\\
f(x_1)+2\sum_{\ell=1}^{r}n_\ell+2r+1\quad & \ r=s,\ldots,t-1\\
2e-3\quad & \ r=t.
\end{cases}$$

\vspace{1cm}
\noindent
In the following we give an explicit construction of a $2$-divisible
$\alpha$-labeling of $H(2t,1)$ for $1<t$ odd, whose existence has been proved in
Theorem \ref{cycle}. We have to distinguish two cases according to the congruence class
of $t$ modulo $4$.

\noindent
Case (1): $t\equiv 1 (mod\ 4)$
\small{$$f(x_r)=\begin{cases}
2t+3-2r\ &  r=1,\ldots,\frac{t+3}{4}\\
2t+2-2r\ &  r=\frac{t+7}{4},\ldots,\frac{t+1}{2}\\
2t+1-2r\ &  r=\frac{t+3}{2},\ldots,t
\end{cases}
\quad\quad
f(y_r^1)=\begin{cases}
2t+2-2r\ &  r=1,\ldots,\frac{t-1}{4}\\
2t+1-2r\ &  r=\frac{t+3}{4},\ldots,\frac{t-1}{2}\\
2t-2r\ &  r=\frac{t+1}{2},\ldots,t
\end{cases}$$}
\small{$$f(x_r^1)=\begin{cases}
3t+2\ &  r=1\\
2t-1+2r\ &  r=2,\ldots,\frac{t+1}{2}\\
3t+1\ &  r=\frac{t+3}{2}\\
2t+2r\ &  r=\frac{t+5}{2},\ldots,t
\end{cases}
\quad\quad
f(y_r)=\begin{cases}
2t+2r\ &  r=1,\ldots,\frac{t-1}{2}\\
3t+3\ &  r=\frac{t+1}{2}\\
2t+1+2r\ &  r=\frac{t+3}{2},\ldots,t
\end{cases}$$}

\noindent
Case (2): $t\equiv 3 (mod\ 4)$

\small{$$f(x_r)=\begin{cases}
2t+3-2r\ &  r=1,\ldots,\frac{t+1}{4}\\
2t+2-2r\ &  r=\frac{t+5}{4},\ldots,\frac{t+1}{2}\\
2t+1-2r\ &  r=\frac{t+3}{2},\ldots,t
\end{cases}
\quad\quad
f(y_r^1)=\begin{cases}
2t+2-2r\ &  r=1,\ldots,\frac{t+1}{4}\\
2t+1-2r\ &  r=\frac{t+5}{4},\ldots,\frac{t+1}{2}\\
2t-2r\ &  r=\frac{t+3}{2},\ldots,t
\end{cases}$$}
\small{$$f(x_r^1)=\begin{cases}
3t+1 &  r=1\\
2t-1+2r &  r=2,\ldots,\frac{t+3}{2}\\
2t+2r &  r=\frac{t+5}{2},\ldots,t,\ \textrm{if}\ t>3
\end{cases}
\quad f(y_r)=\begin{cases}
2t+2r\ &  r=1,\ldots,\frac{t-1}{2}\\
3t+3\ &  r=\frac{t+1}{2}\\
2t+1+2r\  &  r=\frac{t+3}{2},\ldots,t
\end{cases}$$}


\end{document}